\documentclass[reqno,11pt,final]{amsart}
\usepackage{style}
\graphicspath{{./images/}}

\title[Numerical methods for a system of coupled Cahn-Hilliard equations]{Numerical methods for \\ a system of coupled Cahn-Hilliard equations}

\author{M. Martini}
\address{Dipartimento di Matematica, Politecnico di Milano, Via E. Bonardi 9, Milano, Italy}
\email{mattia.martini@mail.polimi.it}

\author{G.E. Sodini}
\address{Dipartimento di Matematica, Politecnico di Milano, Via E. Bonardi 9, Milano, Italy}
\email{giacomoenrico.sodini@mail.polimi.it}

\begin{document}
\begin{abstract}
In this work, we study a system of coupled Cahn-Hilliard equations describing the phase separation of a copolymer and a homopolymer blend.
The numerical methods we propose are based on suitable combinations of existing schemes for the single Cahn-Hilliard equation. As a verification for our approach, we present some tests and a detailed description of the numerical solutions' behaviour obtained by varying the values of the parameters.
\end{abstract}
\maketitle
\tableofcontents
\thispagestyle{empty}

% keywords can be removed
%\keywords{Coupled Cahn-Hilliard equations \and Finite Element Method}

\section{Introduction}
The Cahn-Hilliard equation is a fourth order nonlinear PDE introduced in \cite{cahn, ch2, ch} to model the phase separation in binary systems. This equation and its variants present many practical applications, see for example \cite{vongole, saturno, tumori}.
In the last decades, many methods and numerical techniques to treat this kind of equations have been developed. See \cite{dong} and \cite{teo} for a review of the numerical and theoretical results respectively.\\
The system we study describes the spontaneous separation that occurs in a copolymer and a homopolymer blend. A copolymer consists of two monomers, while a homopolymer is made of a single component. 
The diffuse interface model we consider avoids the explicit treatment of the sharp interface condition between the copolymer and the homopolymer.
This model has been accurately investigated in a three-dimensional setting in \cite{grass}, using the numerical methods developed in \cite{eyre} and \cite{numerical}.\\
In this work we study the two-dimensional, using a different numerical technique, namely Finite Element Method (FEM). In particular, we adapt some basic techniques used for the single Cahn-Hilliard equation and suitably couple them in order to deal with the fully coupled system in an efficient and reliable way. Finally we describe the behavior of the numerical solutions obtained by varying the characteristic parameters of the system. This is done in order to check that the results we obtain reflect the physical properties of the system and performing a stability analysis.\\
The paper is structured as follows: in \autoref{modello} we introduce the mathematical model used to describe the physical phenomenon; in \autoref{numerica} we present the numerical methods for the single equation that will constitute the building blocks for our algorithms. In \autoref{parameters} we describe the behavior of the numerical solutions obtained by varying the parameters of the system. At the end, in \autoref{concl}, we draw some conclusions and we discuss some possible further developments.

\section{The mathematical model}
\label{modello}
In this section we give a brief description of the model whose detailed derivation can be found in \cite{appendice}. The system we consider consists of a copolymer and a homopolymer blend. In particular we have two phase separations: macrophase separation and microphase separation. The first one takes place between the copolymer and the homopolymer and it is represented by the order parameter $u$, which takes values in $[-1,1]$. The ending points of this interval correspond to a homopolymer rich domain $(-1)$ and a copolymer rich domain $(+1)$. The second one is described by the order parameter $v$ and takes place inside the copolymer, between its two component, say A and B. The function $v$ also takes values in [-1,1] with ending points corresponding to A and B rich domains, respectively.\\ 
We focus on the two-dimensional case, mainly due to computational reasons and we study the problem with homogeneous Neumann boundary conditions. Notice that this is a different situation from \cite{grass}, where only periodic boundary conditions are considered. \\
Let $\Omega \subset \mathbb{R}^2$ be a bounded Lipschitz domain and let $\tau_u, \tau_v,  \alpha , \beta, \sigma , \epsilon_u, \epsilon_v, \bar{v}, T \in \mathbb{R}$ with $\epsilon_u, \epsilon_v, \tau_u, \tau_v,$
$T, \sigma >0$. The model is given by the following system of coupled Cahn-Hilliard equations:\\
\begin{equation} \label{sist}
\begin{cases} 
\tau_u u_t = -\Delta \{ \epsilon_u^2 \Delta u + (1-u)(1+u)u -\alpha v -\beta v^2 \}   \quad \quad  & \text{ in } \Omega \times [0, T], \\
\tau_v v_t = -\Delta \{ \epsilon_v^2 \Delta v + (1-v)(1+v)v -\alpha u - 2\beta uv \} -\sigma (v-\overline{v}) \quad \quad  & \text{ in } \Omega \times [0, T],  \\
u( \cdot, 0) = u_0 \quad \quad & \text{ in } \Omega, \\
v( \cdot, 0) = v_0 \quad \quad & \text{ in } \Omega, \\
\partial_{\mathbf{n}}u = \partial_{\mathbf{n}}\{ \epsilon_u^2 \Delta u + (1-u)(1+u)u -\alpha v -\beta v^2 \}=0 & \text{ on } \partial \Omega \times [0, T], \\
\partial_{\mathbf{n}}v = \partial_{\mathbf{n}} \{ \epsilon_v^2 \Delta v + (1-v)(1+v)v -\alpha u - 2\beta uv \} =0  & \text{ on } \partial \Omega \times [0, T].
\end{cases} 
\end{equation}
where $\mathbf{n}$ is the outward normal of $\Omega$ and $u_0, v_0 : \Omega \to \mathbb{R}$ are given initial conditions. \\

Notice that the second equation of the system is a Cahn-Hilliard-Ono equation due to the presence of $\sigma(v-\bar{v})$. The time constants $\tau_u, \tau_v$ control the speed of the evolution of $u$ and $v$ and, specifically, to smaller coefficients correspond a faster evolution. The parameters $\epsilon_u, \epsilon_v$  are proportional to the thickness of the propagation fronts of each component and they control the size of the interface between the macrophases and the microphases.
The parameter $\sigma$ is related to the bonding between A and B in the copolymer and it controls the nonlocal interactions in the second equation. More details on the nonlocal effects in phase separation problems can be found in \cite{otha}. From a practical point of view, values of $\sigma$ different from zero prevent the copolymer from forming large macroscopic domains. As $\sigma$ grows, finer structures and different morphologies arise.
The average value of $v$ in $\Omega$ is represented by $\bar{v}$ and $T$ is the time horizon.
The first coupling parameter $\alpha$ controls the interaction between the confined copolymer and the confining surface. If $\alpha=0$, the configuration is symmetric and thus $u$ has the same preference for any value of $v$. If $\alpha\neq0$, the configuration symmetry is broken, meaning that $u$ has preference for one specific block of the copolymer. The coupling parameter $\beta$ affects the free energy depending of the value of $u$, as $v^2>0$.\\
As explained in \cite{grass}, the state dynamics of these two coupled systems evolves as a gradient flow, up to the reaction term. More precisely, the state variables minimize the value of the following energy functional:
$$ F \equiv F_{\epsilon_u, \epsilon_v, \sigma}(u,v) = \int_{\Omega} \biggl \{ \frac{\epsilon_u^2}{2} | \nabla u|^2 + \frac{\epsilon_v^2}{2} | \nabla v|^2 + W(u,v) + \frac{\sigma}{2} | (-\Delta)^{-1/2}(v-\overline{v})|^2 \biggr \}, $$
where
$$ W(u,v) = \frac{ (u^2-1)^2}{4}+\frac{ (v^2-1)^2}{4} + \alpha uv + \beta uv^2. $$

\section{Numerical approximations}
\label{numerica}
Before presenting our numerical methods for the approximate solution of \eqref{sist}, we revise some known techniques for the single Cahn-Hilliard equation. To this aim, let us set:
\begin{itemize}
\item $ w_u := \epsilon_u^2 \Delta u + (1-u)(1+u)u -\alpha v -\beta v^2$
\item $ w_v := \epsilon_v^2 \Delta v + (1-v)(1+v)v -\alpha u - 2\beta uv$ 
\end{itemize}
and $\phi(x)=(1-x^2)x$. Notice that $w_u$ and $w_v$ are the chemical potentials of the Cahn-Hilliard equations.\\

Then \eqref{sist} becomes
\begin{equation} \label{sistq}
\begin{cases} 
\tau_u u_t = -\Delta w_u \quad \quad  & \text{ in } \Omega \times [0, T], \\
w_u = \epsilon_u^2 \Delta u + \phi(u) -\alpha v -\beta v^2 \quad \quad  & \text{ in } \Omega \times [0, T],  \\
\tau_v v_t = -\Delta w_v -\sigma (v-\overline{v}) \quad \quad  & \text{ in } \Omega \times [0, T],  \\
w_v = \epsilon_v^2 \Delta v + \phi(v) -\alpha u - 2\beta uv  \quad \quad  & \text{ in } \Omega \times [0, T], \\
u( \cdot, 0) = u_0 \quad \quad & \text{ in } \Omega, \\
v( \cdot, 0) = v_0 \quad \quad & \text{ in } \Omega, \\
\partial_{\mathbf{n}}u = \partial_{\mathbf{n}}v = \partial_{\mathbf{n}}w_u = \partial_{\mathbf{n}}w_v =0 & \text{ on } \partial \Omega \times [0, T]. 
\end{cases} 
\end{equation}

\subsection{Semi-discretization in time} \quad \\
In this section we introduce the semi-discretization in time. Let us fix $n>0$, $\Delta t = T/n $ and consider a discretization of $[0,T]$ given by $ \{0=t_0, t_1, \dots, t_{n}, t_{n+1} = T \}$, where $t_k = k\Delta t$ and let $u(t_k) \cong u^{(k)}$ and the same for $v, w_u, w_v$. Consider the following approximation for the time derivative of $u$
\[ u_t \approx \frac{ u^{(k+1)}-u^{(k)} }{\Delta t }, \]
and the same for $v$. By multiplying \eqref{sistq} by the test functions $\varphi, \eta, \psi, \zeta$ and integrating in space over $\Omega$, we obtain: for each $k = 0, \dots, n$ find $(u^{(k)}, v^{(k)}, w_u^{(k)}, w_v^{(k)}  ) \in [H^1(\Omega)]^4$  s.t.
\begin{equation} \label{integral}
\begin{dcases} 
 \int_{\Omega} \tau_u \frac{ u^{(k+1)}-u^{(k)}}{\Delta t} \varphi = \int_{\Omega} \nabla w_u^{(k+1)} \cdot \nabla \varphi, \\
\int_{\Omega} w_u^{(k+1)} \eta = \int_{\Omega} \{ -\epsilon_u^2 \nabla u^{(k+1)} \cdot \nabla \eta + \phi(u^{(k+1)}) \eta -\alpha v^{(k+1)}\eta -\beta (v^{(k+1)})^2 \eta \}, \\
\int_{\Omega} \tau_v \frac{ v^{(k+1)}-v^{(k)}}{\Delta t} \psi = \int_{\Omega} \{ \nabla w^{(k+1)}_v \cdot \nabla \psi -\sigma (v^{(k+1)}-\overline{v})\psi \}, \\
\int_{\Omega} w_v^{(k+1)} \zeta = \int_{\Omega} \{ -\epsilon_v^2 \nabla v^{(k+1)} \cdot \nabla \zeta + \phi(v^{(k+1)}) \zeta -\alpha u^{(k+1)}\zeta -2\beta u^{(k+1)} v^{(k+1)} \zeta \}, \\
u( \cdot, 0) = u_0, \\
v( \cdot, 0) = v_0.
\end{dcases}
\end{equation}
for all $\varphi, \eta, \psi, \zeta  \in H^1(\Omega)$. 
\subsection{Full discretization} \quad \\
To introduce the space discretization of \eqref{integral} based on FEM, we project the semi-discretized problem on a finite dimensional Hilbert space $V_h \subset H^1(\Omega)$. Then the Galerkin approximation of \eqref{integral}, is: for each $k = 0, \dots, n$ find $(u_h^{(k)}, v_h^{(k)}, w_{u,h}^{(k)}, w_{v,h}^{(k)}  ) \in V_h^4$  s.t.
\begin{equation} \label{fulld}
\begin{dcases} 
 \int_{\Omega} \tau_u \frac{ u_h^{(k+1)}-u_h^{(k)}}{\Delta t} \varphi_h = \int_{\Omega} \nabla w_{u,h}^{(k+1)} \cdot \nabla \varphi_h, \\
\int_{\Omega} w_{u,h}^{(k+1)} \eta_h = \int_{\Omega} \{ -\epsilon_u^2 \nabla u_h^{(k+1)} \cdot \nabla \eta_h + \phi(u_h^{(k+1)}) \eta_h -\alpha v_h^{(k+1)}\eta_h -\beta (v_h^{(k+1)})^2 \eta_h \}, \\
\int_{\Omega} \tau_v \frac{ v_h^{(k+1)}-v_h^{(k)}}{\Delta t} \psi_h = \int_{\Omega} \{ \nabla w^{(k+1)}_{v,h} \cdot \nabla \psi_h -\sigma (v_h^{(k+1)}-\overline{v})\psi_h \}, \\
\int_{\Omega} w_{v,h}^{(k+1)} \zeta_h = \int_{\Omega} \{ -\epsilon_v^2 \nabla v_h^{(k+1)} \cdot \nabla \zeta_h + \phi(v_h^{(k+1)}) \zeta_h -\alpha u_h^{(k+1)}\zeta_h -2\beta u_h^{(k+1)} v_h^{(k+1)} \zeta_h \}, \\
u_h( \cdot, 0) = u_{0,h}, \\
v_h( \cdot, 0) = v_{0,h}.
\end{dcases}
\end{equation}
for all $\varphi_h, \eta_h, \psi_h, \zeta_h  \in V_h$, where $u_h^{(k)}, v_h^{(k)}, w_{u,h}^{(k)}, w_{v,h}^{(k)}, u_{0,h}, v_{0,h}$ are the projection of $u^{(k)}, v^{(k)}, w_u^{(k)}, w_v^{(k)}, u_0, v_0$ on $V_h$.  \\

In order to construct such a space, we consider a triangular conforming mesh $\mathcal{T}_h(\Omega)$, where $h=\underset{K\in\mathcal{T}_h}{\max} \diam(K)$. Then a family of possible choices for $V_h$ is
\[
P_r=\{v_h\in C^0(\bar{\Omega})\colon v_h|_K\in\mathbb{P}_r(K),\forall K\in\mathcal{T}_h(\Omega)\}\quad r=1,2,\dots
\]
where $\mathbb{P}_r(K)$ is the space of $r$ degree polynomials defined over $K$. For a more detailed theoretical framework see, for example, Chapter 4 in \cite{quart}.\\

\subsection{Linearization} \label{linear} \quad \\
In this section, we introduce suitable linear approximations for the nonlinear terms $\phi( u^{(k+1)}) $ and $\phi( v^{(k+1)})$ and for the coupling terms $(v^{(k+1)})^2 $ and $u^{(k+1)} v^{(k+1)}$.
\quad \\
For the last two terms we use the following approximations:
\begin{itemize}
\item $(v^{(k+1)})^2  \approx v^{(k)}v^{(k+1)}$,
\item $u^{(k+1)} v^{(k+1)} \approx u^{(k)} v^{(k+1)}$.
\end{itemize}
The first one is the simplest possible choice and seems to be natural. The second choice is motivated by the fact that we expect $u$ to evolve faster then $v$, as pointed out in \cite{grass}. Therefore $u$ will reach a stationary condition before $v$. Then we can reasonably think that the difference $|u^{(k+1)}-u^{(k)}|$ will be smaller than $|v^{(k+1)}-v^{(k)}|$ after the initial macrophase separation; thus, with this strategy, we will achieve a better numerical approximation of the physical phenomenon. \\
\quad \\
\quad \\
Then the key point is to find suitable approximations for the potentials $\phi( u^{(k+1)})$ and $\phi( v^{(k+1)})$. To this aim, we propose the following strategies: 
\begin{enumerate}
\item Optimal dissipation method (OD2) both for $u$ and $v$, i.e.
\[ \phi(u^{(k+1)}) = -\frac{3}{2}(u^{(k)})^2u^{(k+1)}+\frac{1}{2}(u^{(k)})^3 + \frac{u^{(k+1)}+u^{(k)}}{2}, \]
cf. \cite{od2}. Here the authors develop the so-called optimal dissipation approach. In particular a second order in time linear approximation of the potential term $F(u) = 1/4 (\phi^2-1)^2$ is derived by the following Hermite quadrature formula
$$\int_a^b g(x)dx = (b-a)g(a)+\frac{1}{2}(b-a)^2g'(a) + C(b-a)^3 g''(\xi). $$
Using this approximation in a single Cahn-Hilliard equation we get a second order linear scheme.
\item Wu-Van Zwieten-Van der Zee's method (WVV) both for $u$ and $v$:
\[    \phi(u^{(k+1)}) =
\begin{cases}  
2u^{(k+1)} + u^{(k+1)} - u^{(k)} -2 \quad & \text{ if } u^{(k)} < -1, \\
\quad \\
\begin{aligned}
 & 2u^{(k+1)} + u^{(k+1)} - u^{(k)} + 3u^{(k)}-(u^{(k)})^3 \\
 & + \frac{1}{2}(u^{(k+1)}-u^{(k)})(3-3(u^{(k)})^2) \end{aligned} \quad & \text{ if } u^{(k)} \in [-1,1], \\
 \quad \\
2u^{(k+1)} + u^{(k+1)} - u^{(k)} +2  \quad & \text{ if } u^{(k)} > 1.
\end{cases}
\]
which is discussed  in \cite{wvv}. Here the authors consider a diffuse-interface tumor-growth system consisting of a reactive Cahn-Hilliard equation and a reaction-diffusion equation. The schemes are of the Crank-Nicolson type with a convex-concave splitting of the free-energy.
The potential term is approximated employing an implicit Taylor expansion of the convex part and an explicit one of the non-convex part. The splitting considered for the potential term is
$$ F(u) = \begin{dcases}  \biggl ( u^2 + \frac{1}{4} \biggr ) - \biggl ( -2u -\frac{3}{4} \biggr ) \quad & \text{ if } u < -1,\\
\biggl ( u^2 + \frac{1}{4} \bigr ) - \biggl ( \frac{3}{2}u^2 -\frac{1}{4} u^4 \biggr ) \quad & \text{ if } u \in [-1,1], \\
\biggl ( u^2 + \frac{1}{4} \biggr ) - \biggl (2u -\frac{3}{4} \biggr ) & \text{ if } u > 1. \\
\end{dcases}
$$
Moreover, this method stabilizes the system by modifying the second and the fourth equations of \eqref{integral} as follows:
\begin{align*}
\int_{\Omega} w_u^{(k+1)} \psi  &= \int_{\Omega} \{ -\frac{\epsilon^2}{2} \nabla u^{(k+1)} \cdot \nabla \psi + \phi(u^{(k+1)})\psi \} \\
& + \biggl [ \int_{\Omega} \{ -\frac{\epsilon^2}{2} \nabla u^{(k)} \cdot \nabla \psi  - \alpha \nabla (u^{(k+1)} -u^{(k)}) \cdot \nabla \psi \} \biggr ], \\
\end{align*}
where $ \alpha>0 $ is a stabilization parameter. The expression in squared brackets is completely new and in addition there is a weighted contribution of $ \nabla u^{(k+1)}$ which is multiplied by $\frac{1}{2}$. Clearly the same applies to $v$.
According to \cite{od2}, this method is second order accurate both in space and time in the single Cahn-Hilliard equation case.
\item Eyre's method (EY), both for $u$ and $v$:
\[ \phi(u^{(k+1)}) = -2 u^{(k+1)} - (u^{(k)})^3 +3u^{(k)}, \]
The key point in Eyre's work \cite{eyr} is to approximate the potential term by introducing a positive phobic numerical dissipation in the discrete energy law, to ensure the unconditional energy-stability of the scheme.
\item Linear splitting method (LS) both for $u$ and $v$; i.e.
\[ \phi(u^{(k+1)})= (1-u^{(k)})(1+u^{(k)})u^{(k+1)} \]
suggested in \cite{grass}.To treat nonlinearities, the cubic term is split into the product of two terms; the first one is quadratic and related to the state of the system at the present time step and the second one is linear and related to the state of the system at the next time step.
\end{enumerate}
\quad \\
\quad \\
\subsection{OD2, WVV, LS and EY methods: a comparison} 
\begin{figure}[h!]
    \centering
    \includegraphics[scale=0.1]{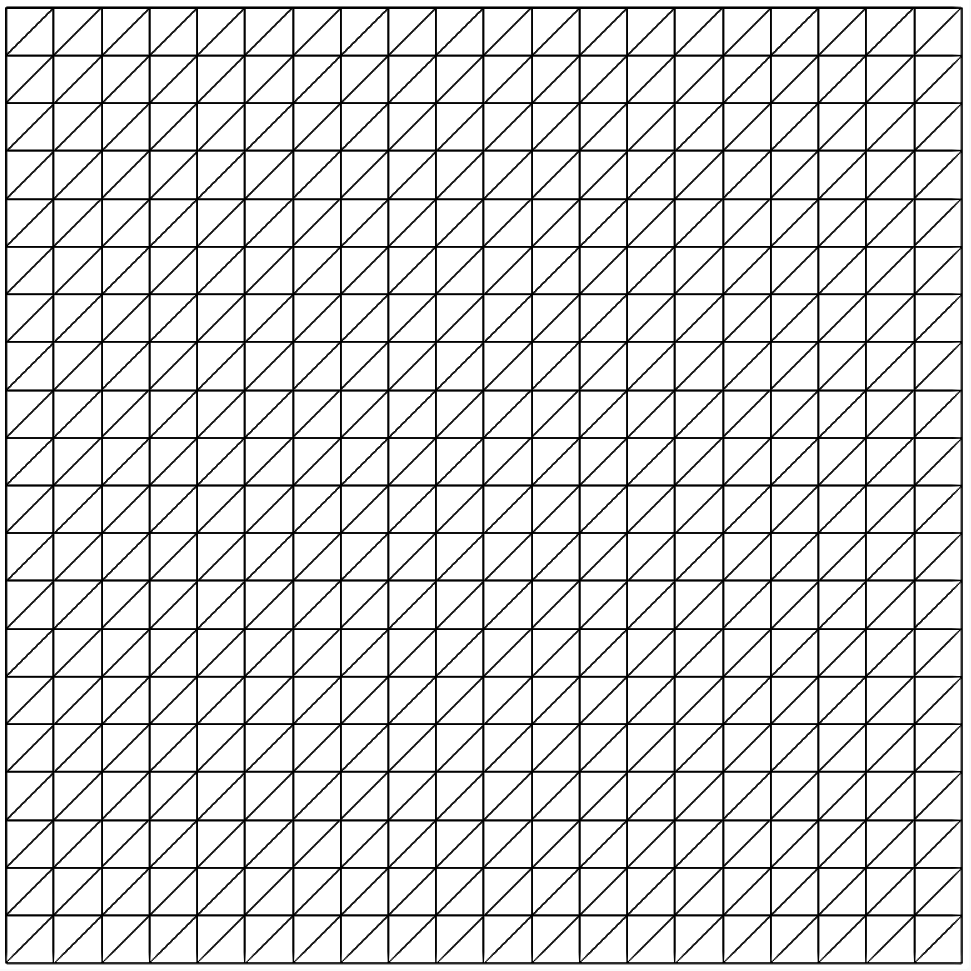}
    \caption{$20\times20$ conforming triangular mesh used for the tests}
    \label{th}
\end{figure}

\begin{figure}[h!]
	\centering
	\includegraphics[height=4.5cm, width=4.5cm]{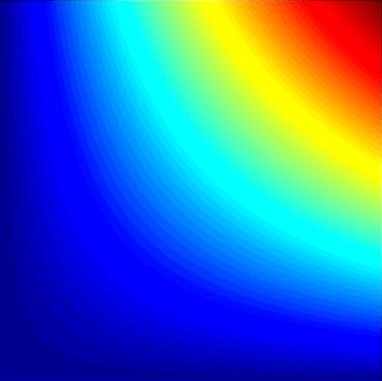} \quad \quad
	\includegraphics[height=4.69cm, width=5.45cm]{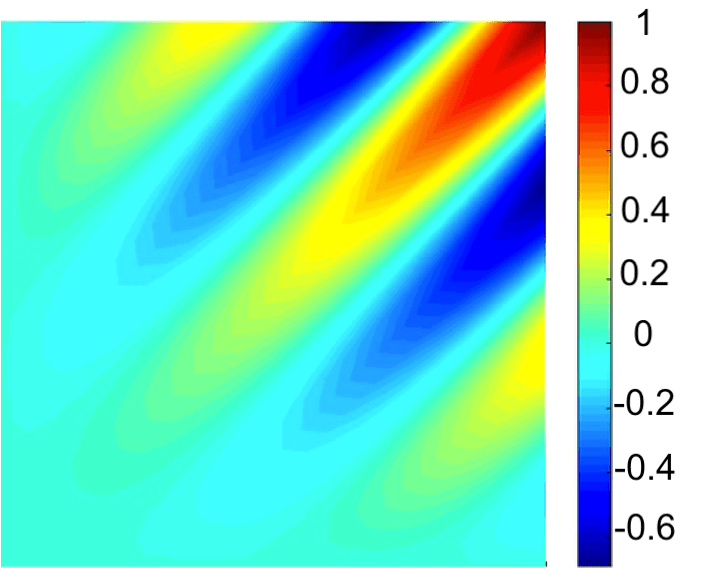}
	\caption{Initial conditions of $u$ (left) and $v$ (right).}
		\label{initcond}
	\end{figure}

\quad \\
We first verify that, when the coupling parameters $\alpha$ and $\beta$ are set equal to $0$, the computed solutions (for $u$ and $v$ taken separately) coincide with the ones of the single equation case. After that, we set nonzero values to the coupling parameters and we analyze the results obtained by the different numerical methods detailed in \autoref{linear}.
In the following test, the time step of each method is adjusted in order to have stable numerical solutions, in particular for OD2 $\Delta t = 5 \cdot10^{-3}$, for LS $\Delta t = 10^{-4}$ and for EY $\Delta t = 10^{-4}$. Notice that the time step used for the three methods are different; in particular, the one for OD2 is larger than the others. The setting used for the simulations is the following: \\ $\Omega = (0,1) \times (0,1)$, $u_0= \sin(xy)$, $\tau_u =1$, $v_0= \cos(10(x-y))xy$, $\tau_v=1$, $\sigma=0.3$, $\overline{v} = \int_{\Omega} v_0 \approx 0.0114559 $, $\epsilon_u = 0.05$, $\epsilon_v=0.05$, $\alpha=0.5$, $\beta=0.8$, $T=10$. In all the following simulations the chosen FEM space is $P_1$ and the mesh $\mathcal{T}_h$ is a $20\times 20$ grid as in \autoref{th}.
Notice that this parameters have been chosen to verify the methods. In \autoref{initcond} the snapshots of the initial conditions are reported. As we can see in \autoref{u3test}, $u$ and $v$ are very similar, both in the shape of the solution and in the numerical values attained. Then we choose to use the OD2 method to perform tests on parameters. Since this method achieves reasonably accurate numerical results even with large temporal steps (in some cases $0.01$), it is easier from a computational point of view to use it for tests on parameters.

	\quad \\
	
\begin{table}[h!]
  \centering
  \begin{tabular}{  m{1cm}  m{5cm} m{5cm}  m{5.3cm} }
    \hline
     Method & \quad \quad \quad  \quad   $t=0.2$ & \quad \quad \quad \quad  $t=2$  &   \quad \quad \quad \quad  $t=10$ \\ \hline
     \quad & \quad & \quad & \quad 
     \\
      \quad OD2
    & 
    \begin{minipage}{.3\textwidth}
      \includegraphics[width=1.9cm, height=1.9cm]{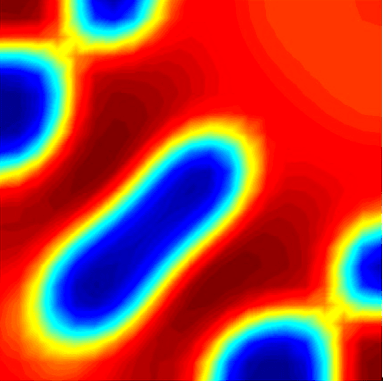}
      \includegraphics[width=1.9cm, height=1.9cm]{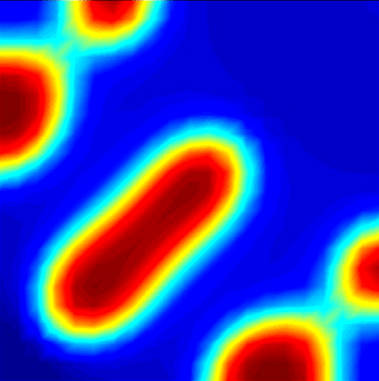}
    \end{minipage}
    &
     \begin{minipage}{.3\textwidth}
       \includegraphics[width=1.9cm, height=1.9cm]{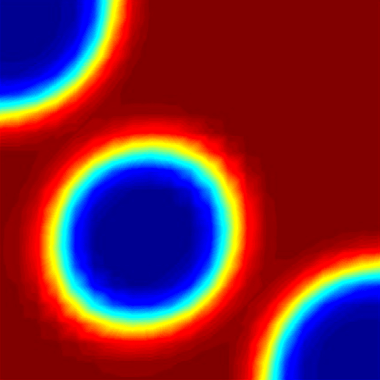}
      \includegraphics[width=1.9cm, height=1.9cm]{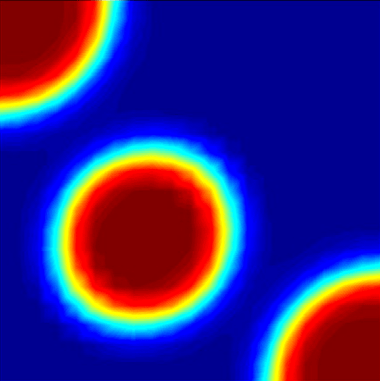}
    \end{minipage}
    &
     \begin{minipage}{.3\textwidth}
      \includegraphics[width=1.9cm, height=1.9cm]{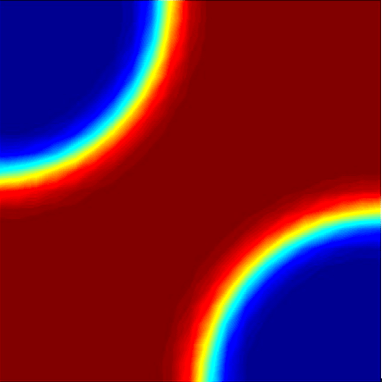}
      \includegraphics[width=1.9cm, height=1.9cm]{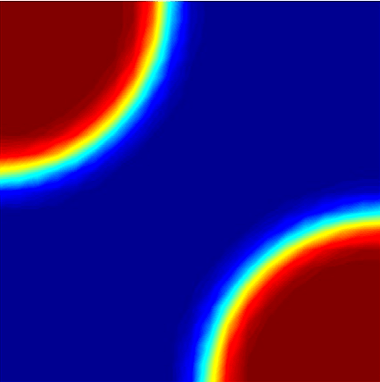}
    \end{minipage} \\
    \quad & \quad & \quad & \quad \\
    \hline
     \quad & \quad & \quad & \quad \\
    \quad LS
    & 
    \begin{minipage}{.3\textwidth}
      \includegraphics[width=1.9cm, height=1.9cm]{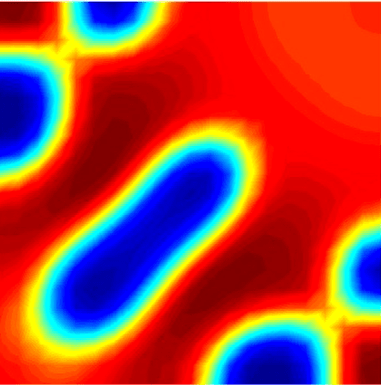}
      \includegraphics[width=1.9cm, height=1.9cm]{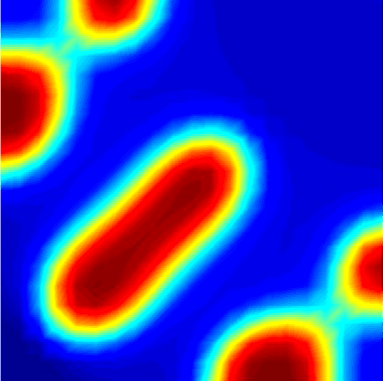}
    \end{minipage}
    &
     \begin{minipage}{.3\textwidth}
     \includegraphics[width=1.9cm, height=1.9cm]{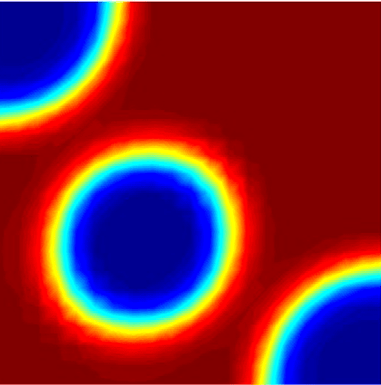}
      \includegraphics[width=1.9cm, height=1.9cm]{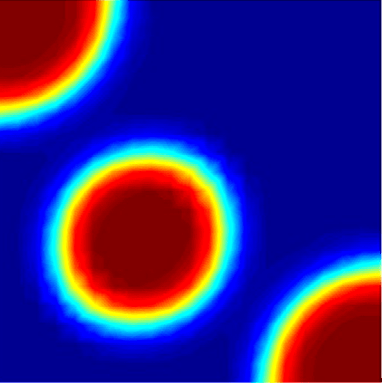}
    \end{minipage}
    &
     \begin{minipage}{.3\textwidth}
      \includegraphics[width=1.9cm, height=1.9cm]{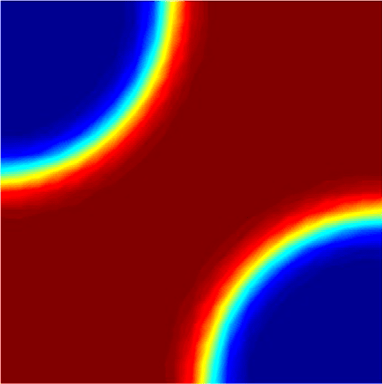}
      \includegraphics[width=1.9cm, height=1.9cm]{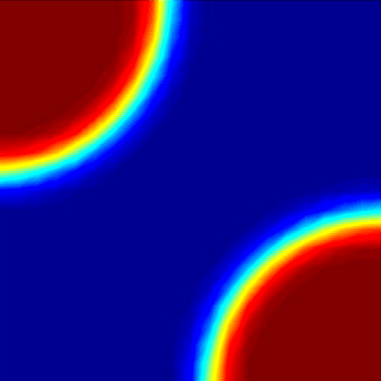}
    \end{minipage} 
    \\ 
    \quad & \quad & \quad & \quad 
    \\
    \hline
     \quad & \quad & \quad & \quad 
    \\
    \quad EY
    & 
    \begin{minipage}{.3\textwidth}
       \includegraphics[width=1.9cm, height=1.9cm]{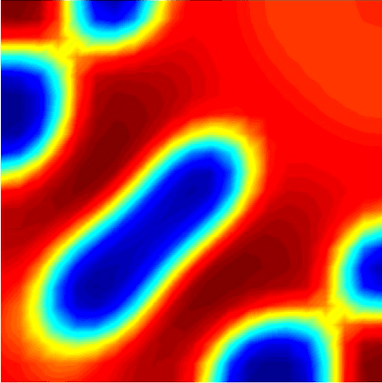}
      \includegraphics[width=1.9cm, height=1.9cm]{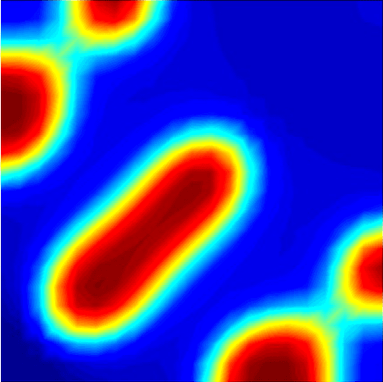}
    \end{minipage}
    &
     \begin{minipage}{.3\textwidth}
        \includegraphics[width=1.9cm, height=1.9cm]{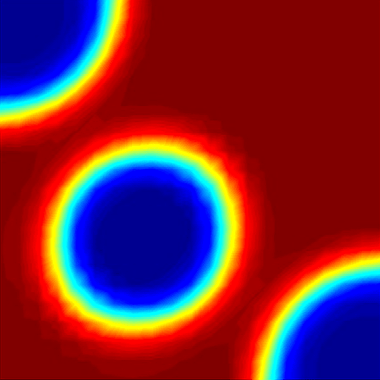}
      \includegraphics[width=1.9cm, height=1.9cm]{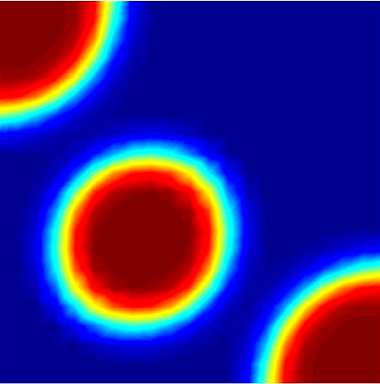}
    \end{minipage}
    &
     \begin{minipage}{.3\textwidth}
       \includegraphics[width=1.9cm, height=1.9cm]{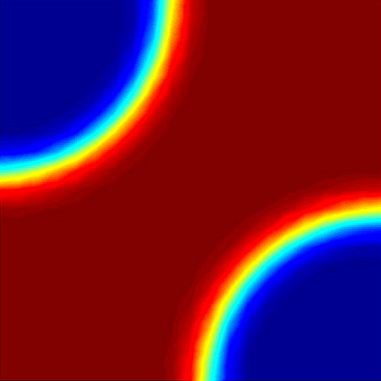}
      \includegraphics[width=1.9cm, height=1.9cm]{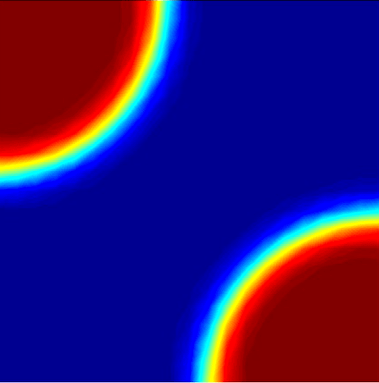}
    \end{minipage} 
    \\
    \quad & \quad & \quad & \quad \\
    
    \hline
        \quad & 
        \quad\quad\quad\quad\quad\quad\quad\quad\quad\quad\quad 
\begin{minipage}{.3\textwidth}
       \includegraphics[trim=0 0 0 -2, scale=0.35]{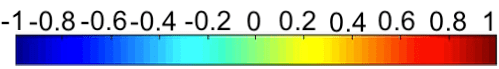}
    \end{minipage}
    &
    \quad 
    &
    \quad \\

  \end{tabular}
  \caption{Comparison of the evolution of $u$ and $v$.} 
  \label{u3test}
\end{table}	

We conclude that
\begin{enumerate}
	
\item OD2 method gives encouraging results: the numerical solutions seem to be physically relevant, even with "large" time steps (order $10^{-2}$). We analyze in further details the output of this method in the forthcoming section;
\item Method WVV does not give reasonable outputs: even with time steps of $10^{-6}$ and high number of nodes, the numerical solutions tend immediately to reach very high values (order $10^6$) and to oscillate very fast. We conclude that this approximation is not suitable for the coupled case, as this behavior is not physically acceptable;
\item Methods LS and EY give accurate results only for sufficiently small time steps (order $10^{-4}$). 
\end{enumerate}
 \section {Fully coupled system}
 \label{parameters}
\begin{figure}[h!] 
\centering
\includegraphics[height=5cm,width=13cm]{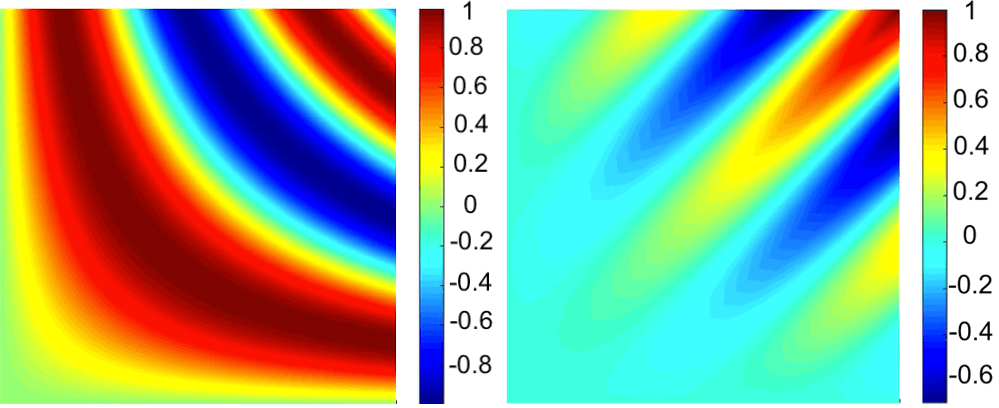}
	\caption{Initial conditions of $u$ (left) and $v$ (right).}
	\label{initcond2}
\end{figure}

In this section, we verify that the numerical solutions we obtained by varying the parameters respect the physical properties of the system. To achieve this purpose, we present some tests, performed with OD2 linearization method. We remark that this method allows larger time steps than EY and LS, speeding up the computations.\\
For each test the general setting is the following: $\Omega= (0,1) \times (0,1)$, $u_0= \sin(10xy)$, $v_0= \cos(10(x-y))xy$, $\overline{v} = \int_{\Omega} v_0 \approx 0.0114559 $, $\Delta t=0.005$, $T=15$ and the FEM space used is $P1$ with a $20\times20$ grid. For each parameter we report in a figure the behavior of the solution at different times and for different values of the parameters. The plots of the initial conditions for all test are reported in \autoref{initcond2} and the values chosen for our tests are inspired by \cite{appendice}. In particular we vary them in the ranges considered by the authors.\\
\quad \\
  \textbf{Test on $\tau_v$}\\
  In this test we vary the value of $\tau_v$, in order to verify that its growth implies a slower evolution of the copolymer part of the system.  As we can see in \autoref{table: tv} 
 the system reaches the same final configuration for both values of the parameter, while in the case $\tau_v=250$ it evolves slower: for $\tau_v = 10$ the part of the copolymer in the upper-right corner has already disappered at time $t = 3$, while for $\tau_v = 250$, it is still present at time t = 6. The same happens for $\tau_u$.
 \begin{table}[h!]
  \centering
  \begin{tabular}{   m{1cm}  m{5cm} m{5cm}  m{5.3cm} }
    \hline
     $\quad \tau_v$ &  \quad \quad \quad  \quad   $t=3$ &   \quad \quad \quad \quad $t=6$  & \quad \quad \quad \quad $t=15$ \\ \hline
     \quad & \quad & \quad & \quad 
     \\
      \quad $10$
    & 
    \begin{minipage}{.3\textwidth}
      \includegraphics[width=1.9cm, height=1.9cm]{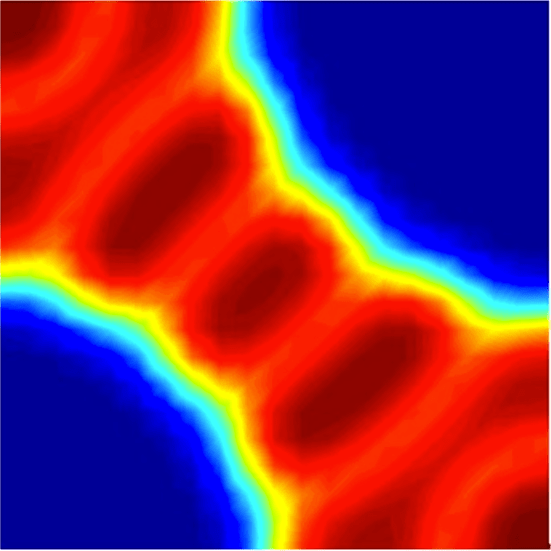}
      \includegraphics[width=1.9cm, height=1.9cm]{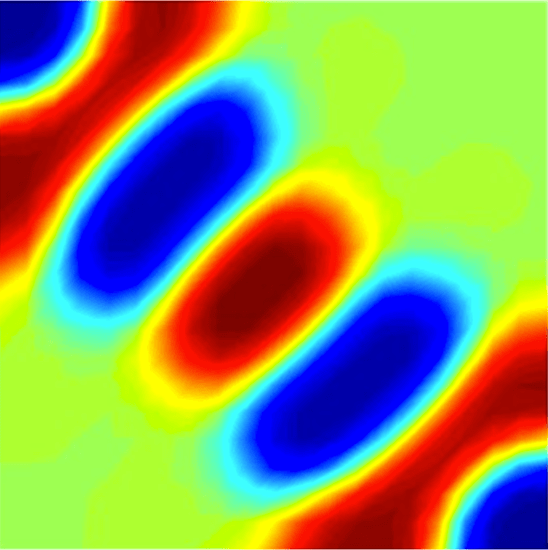}
    \end{minipage}
    &
     \begin{minipage}{.3\textwidth}
       \includegraphics[width=1.9cm, height=1.9cm]{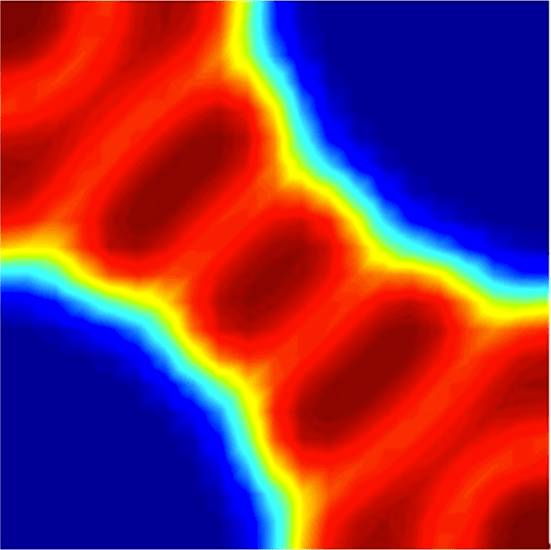}
      \includegraphics[width=1.9cm, height=1.9cm]{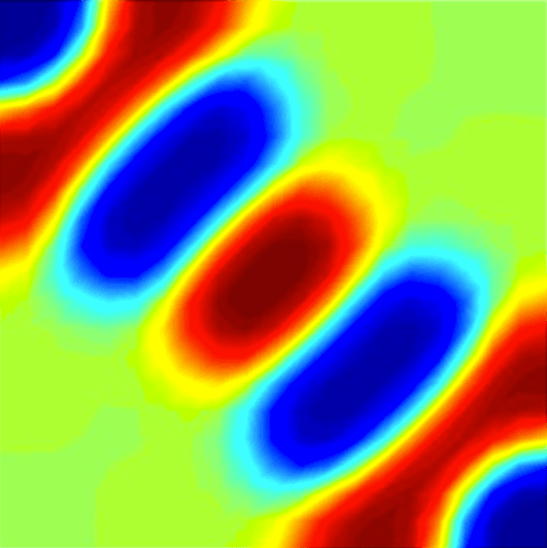}
    \end{minipage}
    &
     \begin{minipage}{.3\textwidth}
             \includegraphics[width=1.9cm, height=1.9cm]{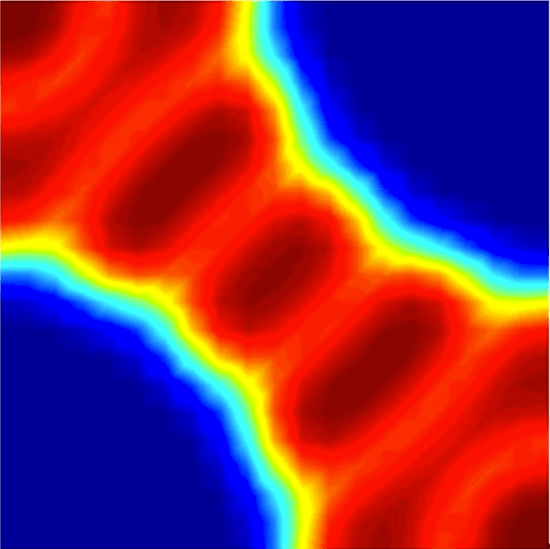}
      \includegraphics[width=1.9cm, height=1.9cm]{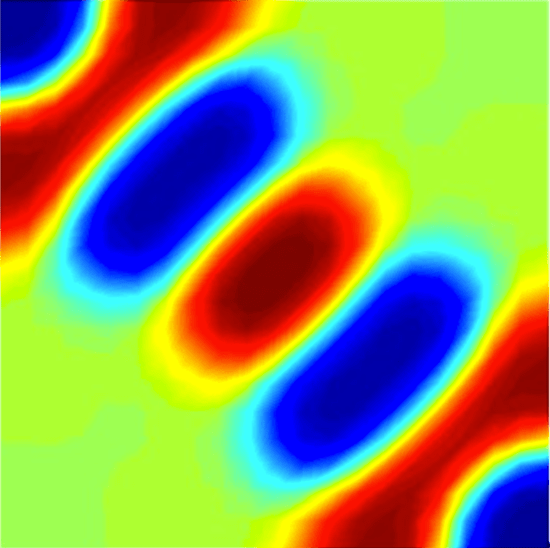}
    \end{minipage}
    \\ 
    \quad & \quad & \quad & \quad 
    \\
   
    \hline
     \quad & \quad & \quad & \quad 
    \\
    \quad $250$
    & 
    \begin{minipage}{.3\textwidth}
            \includegraphics[width=1.9cm, height=1.9cm]{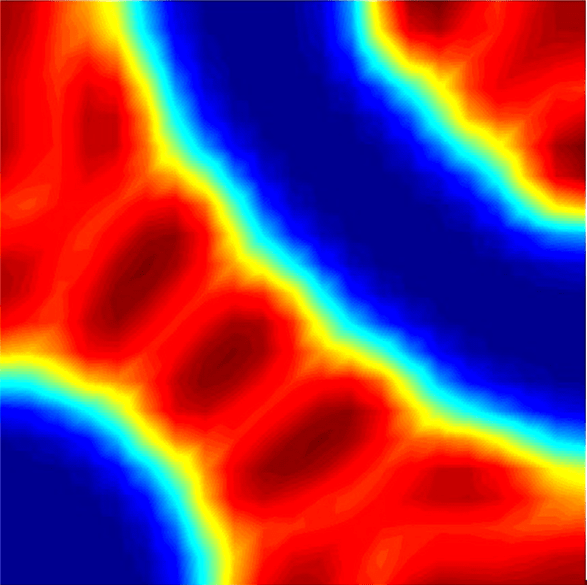}
      \includegraphics[width=1.9cm, height=1.9cm]{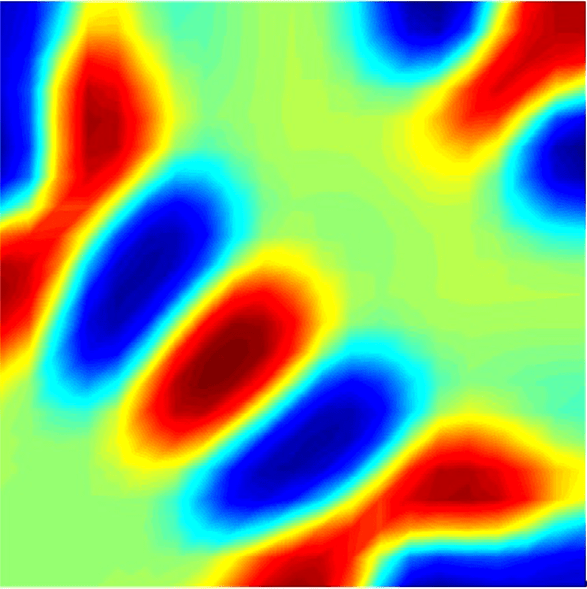}
    \end{minipage}
    &
     \begin{minipage}{.3\textwidth}
     \includegraphics[width=1.9cm, height=1.9cm]{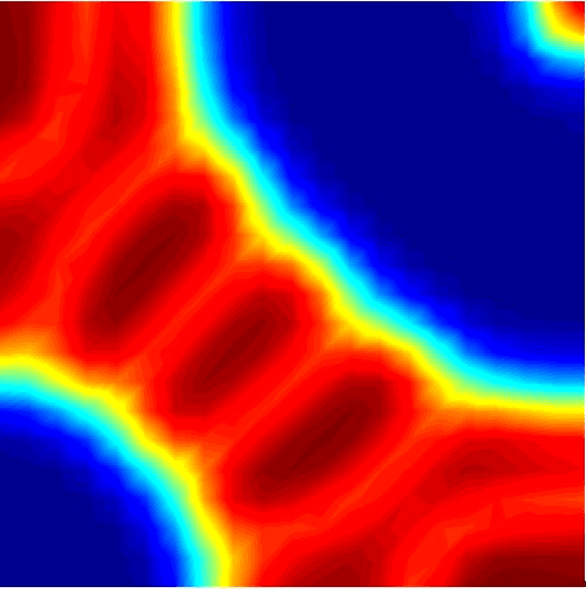}
      \includegraphics[width=1.9cm, height=1.9cm]{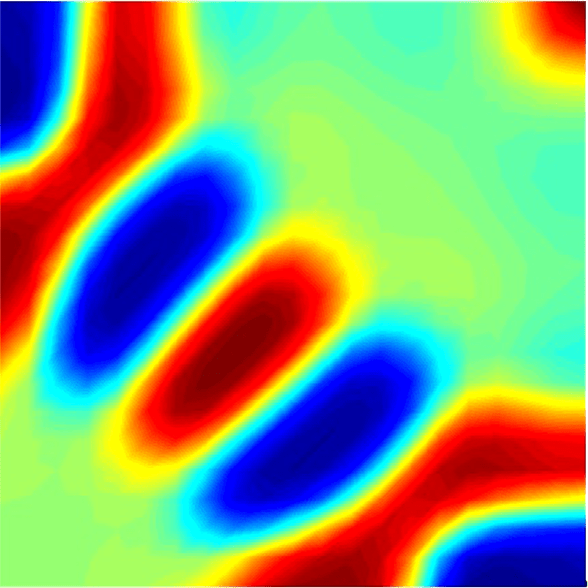}
    \end{minipage}
    &
     \begin{minipage}{.3\textwidth}
      \includegraphics[width=1.9cm, height=1.9cm]{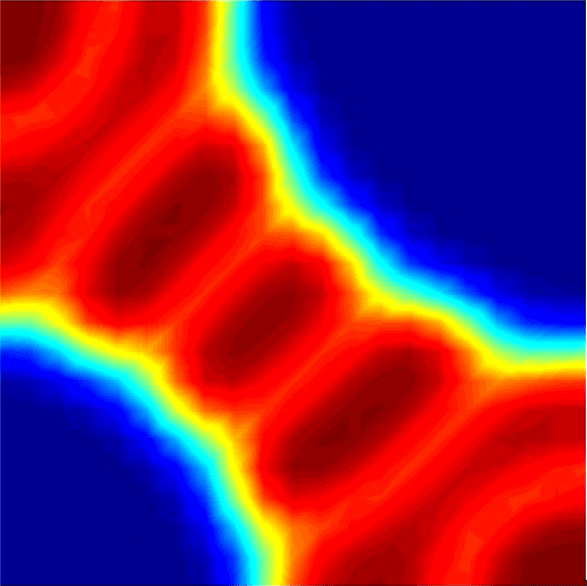}
      \includegraphics[width=1.9cm, height=1.9cm]{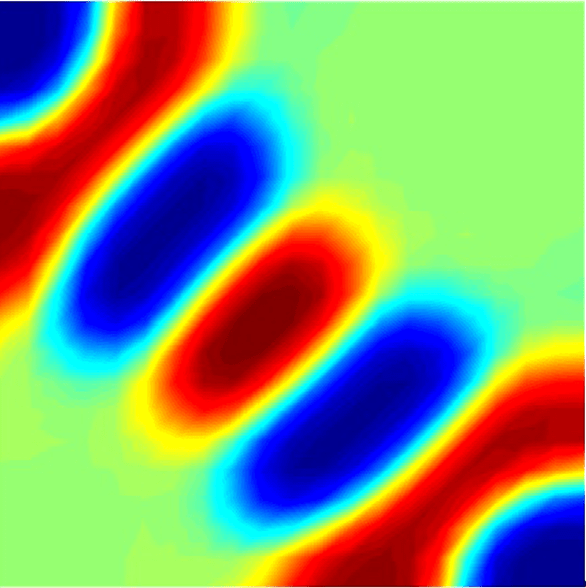}
    \end{minipage} 
    \\
    \quad & \quad & \quad & \quad \\

    \hline 
        \quad &
\quad\quad\quad\quad\quad\quad\quad\quad\quad\quad\quad 
\begin{minipage}{.3\textwidth}
       \includegraphics[trim=0 0 0 -2, scale=0.35]{colorbar.png}
    \end{minipage}
    &
    \quad 
    &
    \quad \\
  \end{tabular}
  \caption{Evolution of $u$ (left) and $v$ (right) with parameters $\epsilon_u=\epsilon_v=0.05$, $\tau_u=1$, $\sigma=100$, $\alpha=0.04$ and $\beta=-0.9$.}
  \label{table: tv}
\end{table}
\pagebreak
\quad \\
\quad \\
\textbf{Test on $\epsilon_v$} \\
 This parameter controls the separation interface thickness between the two components of the copolymer. For bigger values of $\epsilon_v$ the interface between the two copolymers must be thicker so that it prevents the generation of fine patterns. As we can see in \autoref{table: ev}, for small values of $\epsilon_v$, we have more complicated patterns which become coarser, according to the growth of $\epsilon_v$. The same happens for $\epsilon_u$.\\
  \begin{table}[h!]
  \centering
  \begin{tabular}{  m{1cm}  m{5cm} m{5cm}  m{5.3cm} }
    \hline
     $\quad \epsilon_v$ &  \quad \quad \quad  \quad   $t=0.9$ &   \quad \quad \quad \quad  $t=6$  &    \quad \quad \quad \quad  $t=15$ \\ \hline
     \quad & \quad & \quad & \quad 
     \\
      \quad $0.01$
    & 
    \begin{minipage}{.3\textwidth}
      \includegraphics[width=1.9cm, height=1.9cm]{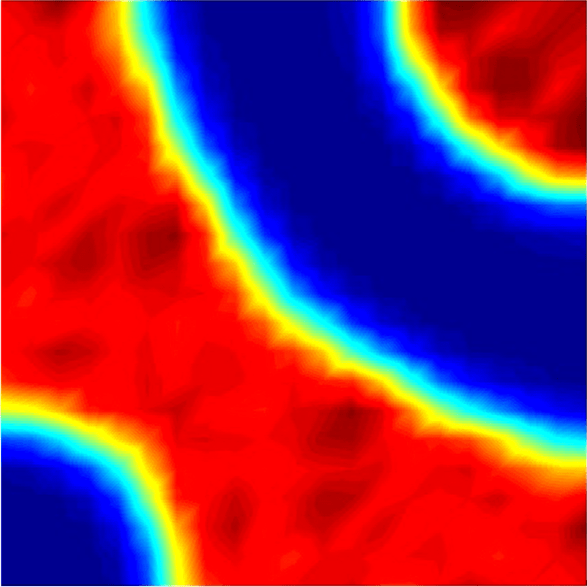}
      \includegraphics[width=1.9cm, height=1.9cm]{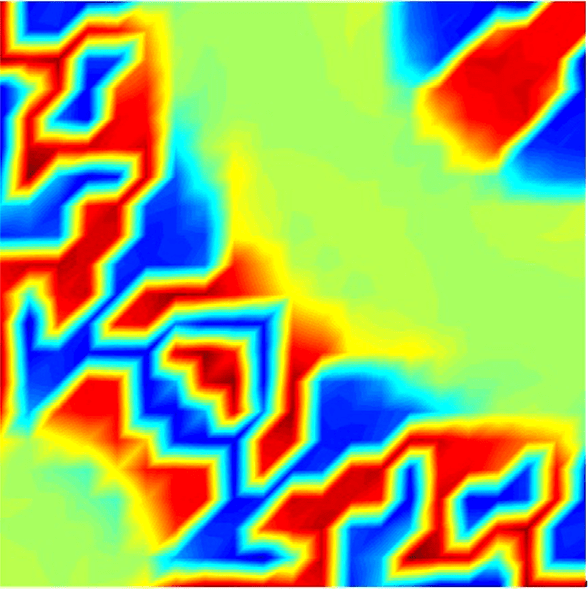}
    \end{minipage}
    &
     \begin{minipage}{.3\textwidth}
       \includegraphics[width=1.9cm, height=1.9cm]{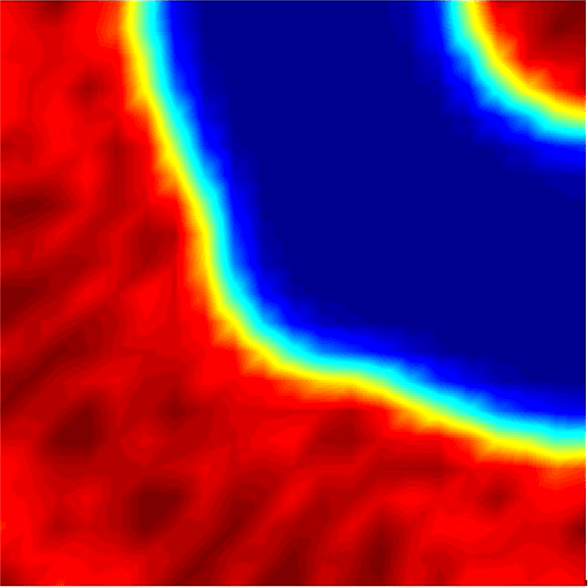}
      \includegraphics[width=1.9cm, height=1.9cm]{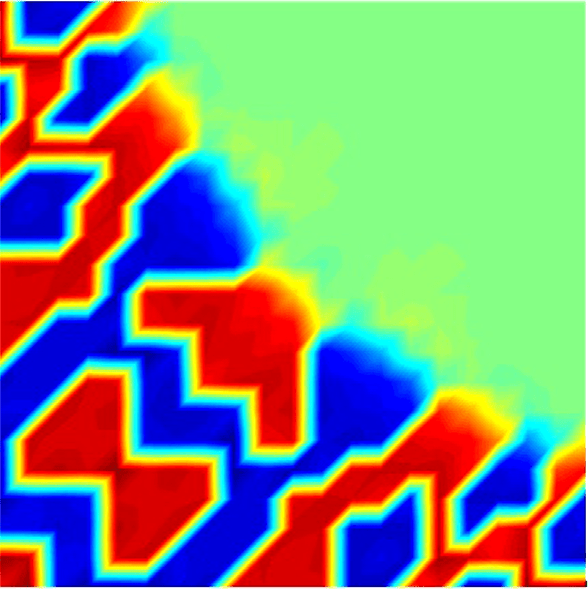}
    \end{minipage}
    &
     \begin{minipage}{.3\textwidth}
      \includegraphics[width=1.9cm, height=1.9cm]{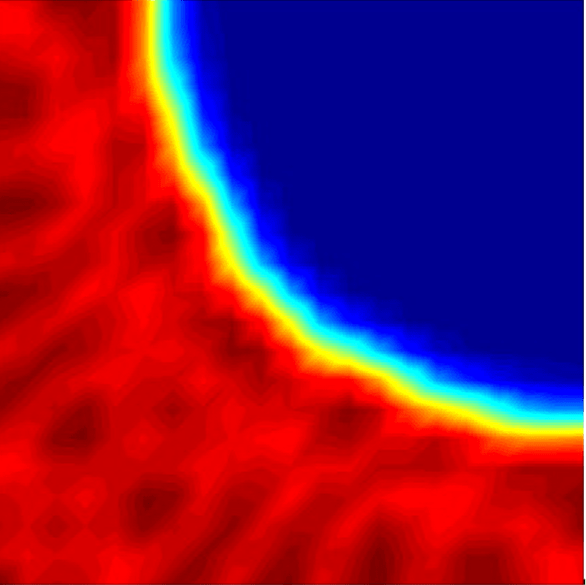}
      \includegraphics[width=1.9cm, height=1.9cm]{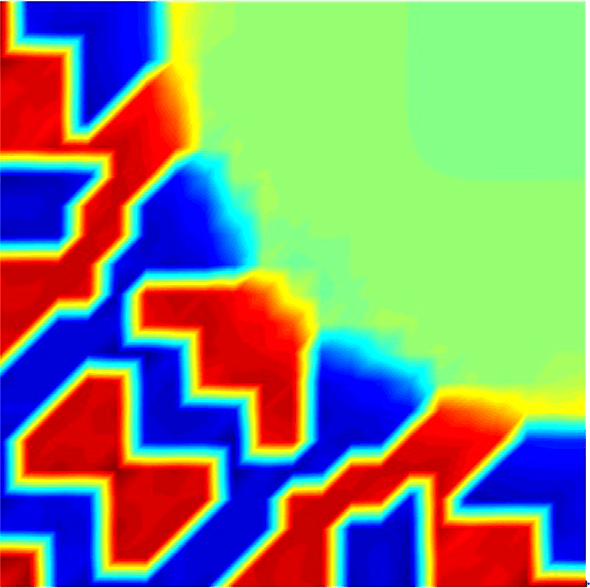}
    \end{minipage} \\
    \quad & \quad & \quad & \quad \\
    \hline
     \quad & \quad & \quad & \quad \\
    \quad $0.03$
    & 
    \begin{minipage}{.3\textwidth}
      \includegraphics[width=1.9cm, height=1.9cm]{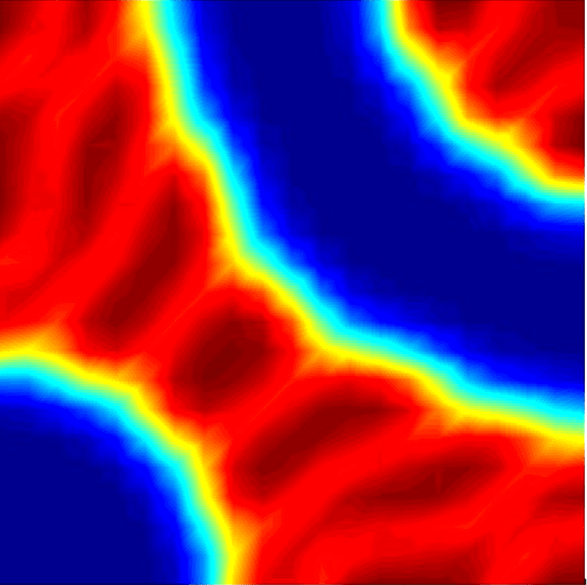}
      \includegraphics[width=1.9cm, height=1.9cm]{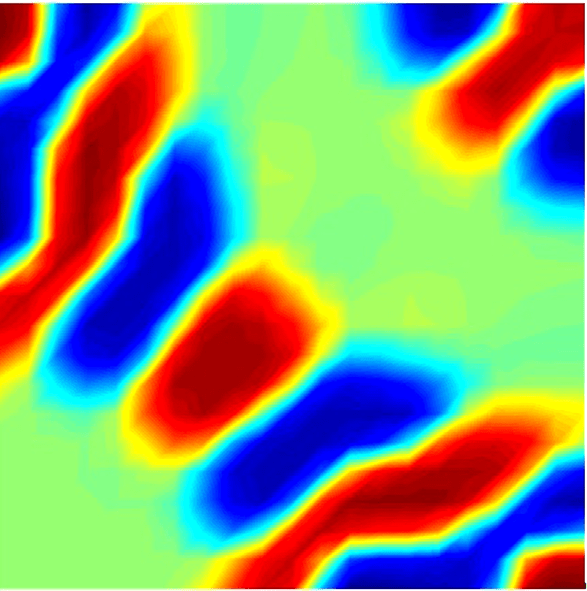}
    \end{minipage}
    &
     \begin{minipage}{.3\textwidth}
     \includegraphics[width=1.9cm, height=1.9cm]{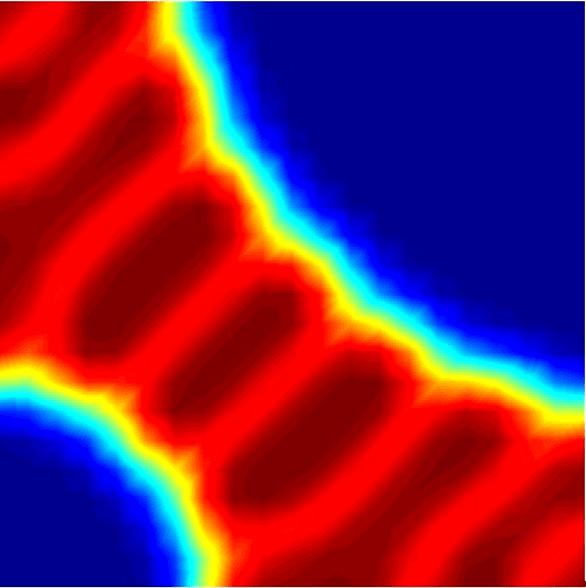}
      \includegraphics[width=1.9cm, height=1.9cm]{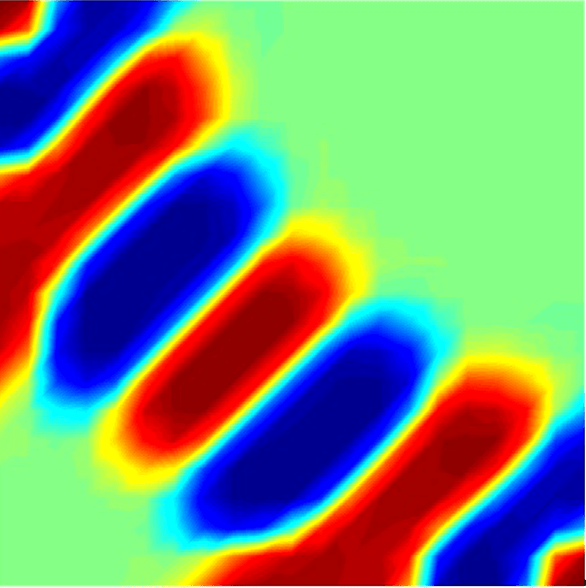}
    \end{minipage}
    &
     \begin{minipage}{.3\textwidth}
      \includegraphics[width=1.9cm, height=1.9cm]{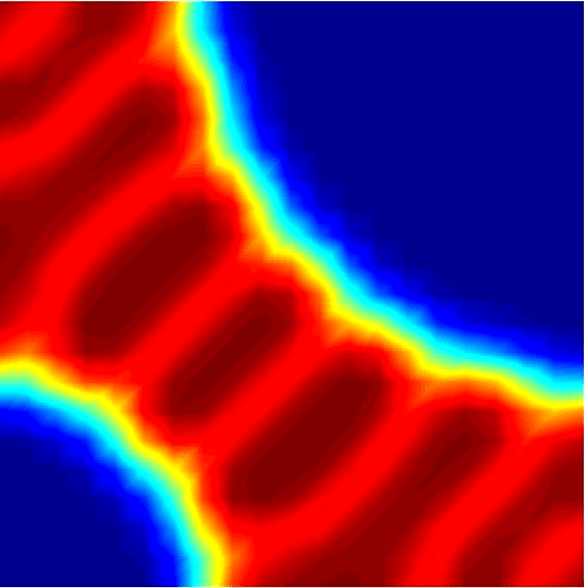}
      \includegraphics[width=1.9cm, height=1.9cm]{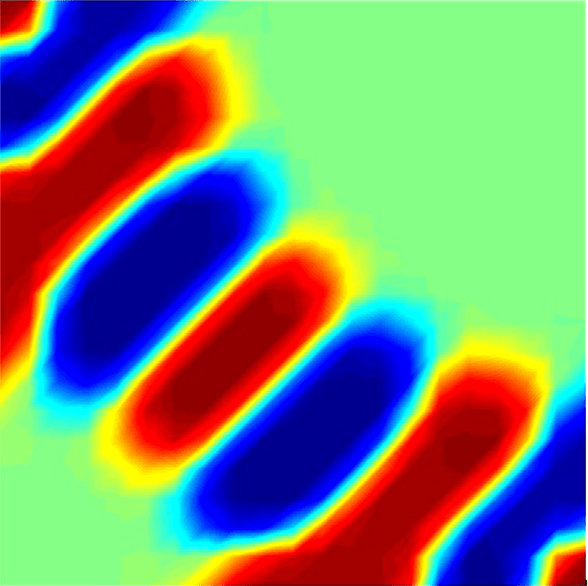}
    \end{minipage} 
    \\ 
    \quad & \quad & \quad & \quad 
    \\
    \hline
     \quad & \quad & \quad & \quad 
    \\
    \quad $0.05$
    & 
    \begin{minipage}{.3\textwidth}
       \includegraphics[width=1.9cm, height=1.9cm]{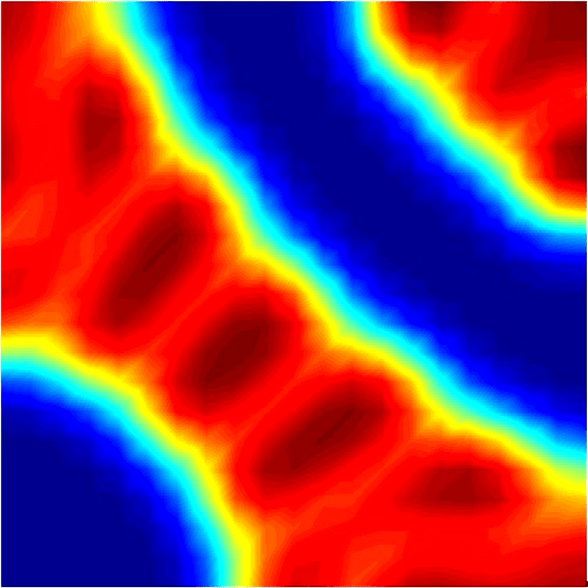}
      \includegraphics[width=1.9cm, height=1.9cm]{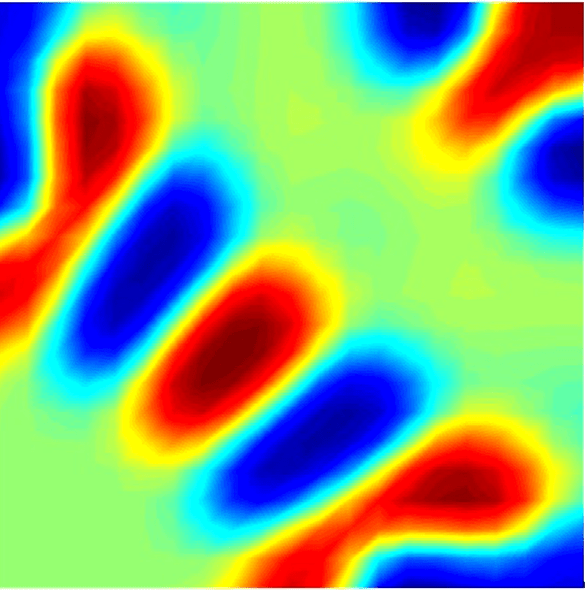}
    \end{minipage}
    &
     \begin{minipage}{.3\textwidth}
        \includegraphics[width=1.9cm, height=1.9cm]{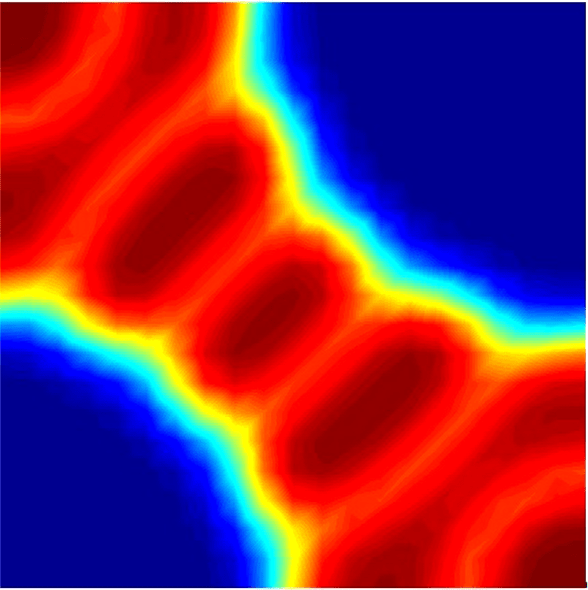}
      \includegraphics[width=1.9cm, height=1.9cm]{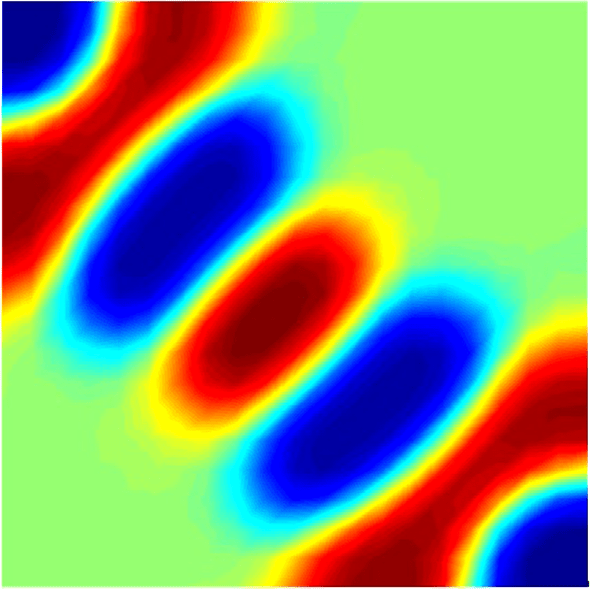}
    \end{minipage}
    &
     \begin{minipage}{.3\textwidth}
       \includegraphics[width=1.9cm, height=1.9cm]{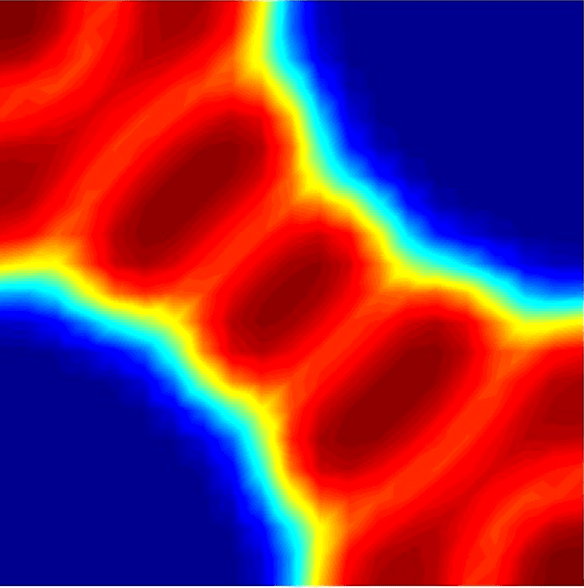}
      \includegraphics[width=1.9cm, height=1.9cm]{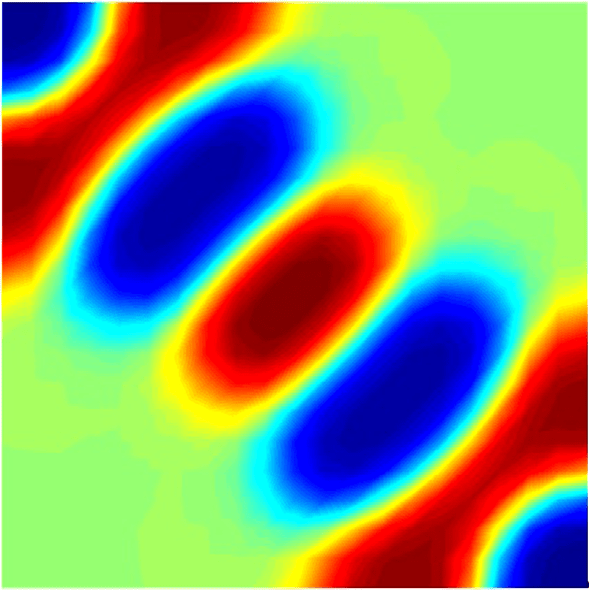}
    \end{minipage} 
    \\
    \quad & \quad & \quad & \quad \\
    
    \hline
        \quad &
\quad\quad\quad\quad\quad\quad\quad\quad\quad\quad\quad 
\begin{minipage}{.3\textwidth}
       \includegraphics[trim=0 0 0 -2, scale=0.35]{colorbar.png}
    \end{minipage}
    &
    \quad 
    &
    \quad \\
  \end{tabular}
  \caption{Evolution of $u$ (left) and $v$ (right) with parameters $\epsilon_u=0.05$, $\tau_u=1$, $\tau_v=100$, $\sigma=60$, $\alpha=0.02$ and $\beta=-0.9$.} 
  \label{table: ev}
\end{table}
\quad \\
\pagebreak
\quad \\
\quad \\
\textbf{Test on $\sigma$}\\
  The parameter $\sigma$ is related to the connectivity between the two components of the copolymer. As we can see in \autoref{table: sigma}, as the value of $\sigma$ increases, the components of the copolymer tend no longer to form big macroareas, but they separate in more and smaller parts.

 \begin{table}[h!]
  \centering
  \begin{tabular}{   m{1cm}  m{5cm} m{5cm}  m{5.3cm} }
    \hline
     $\quad \sigma$ &  \quad \quad \quad  \quad   $t=3$ &   \quad \quad \quad \quad  $t=6$  &    \quad \quad \quad \quad  $t=15$ \\ \hline
     \quad & \quad & \quad & \quad 
     \\
      \quad $0$
    & 
    \begin{minipage}{.3\textwidth}
      \includegraphics[width=1.9cm, height=1.9cm]{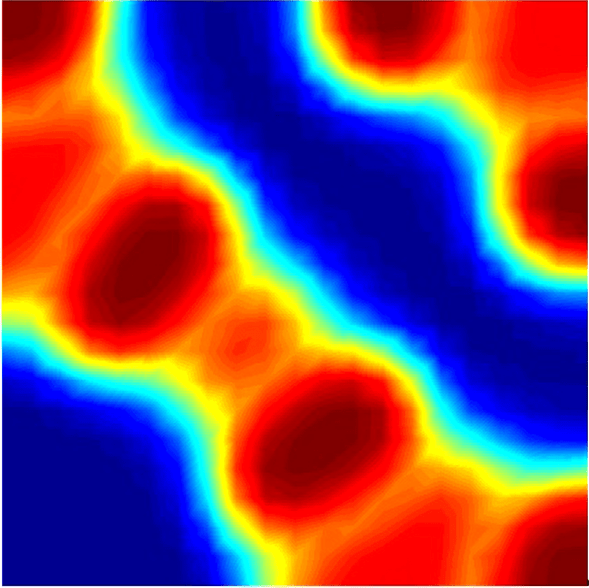}
      \includegraphics[width=1.9cm, height=1.9cm]{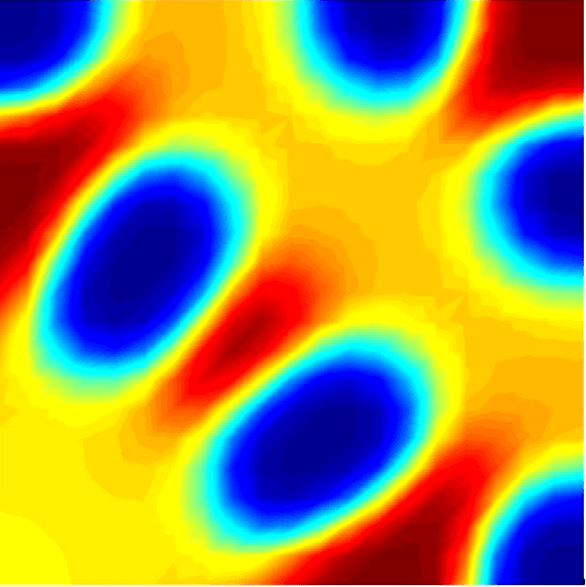}
    \end{minipage}
    &
     \begin{minipage}{.3\textwidth}
       \includegraphics[width=1.9cm, height=1.9cm]{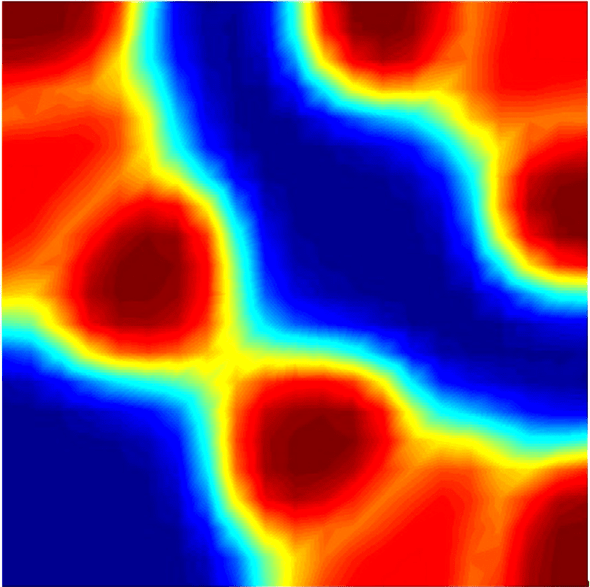}
      \includegraphics[width=1.9cm, height=1.9cm]{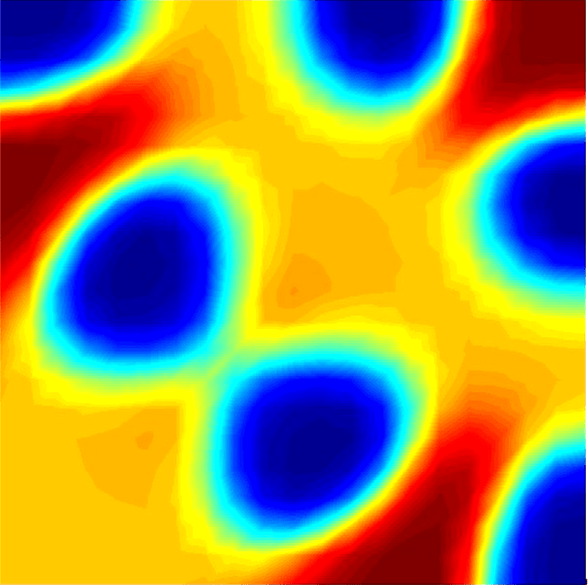}
    \end{minipage}
    &
     \begin{minipage}{.3\textwidth}
        \includegraphics[width=1.9cm, height=1.9cm]{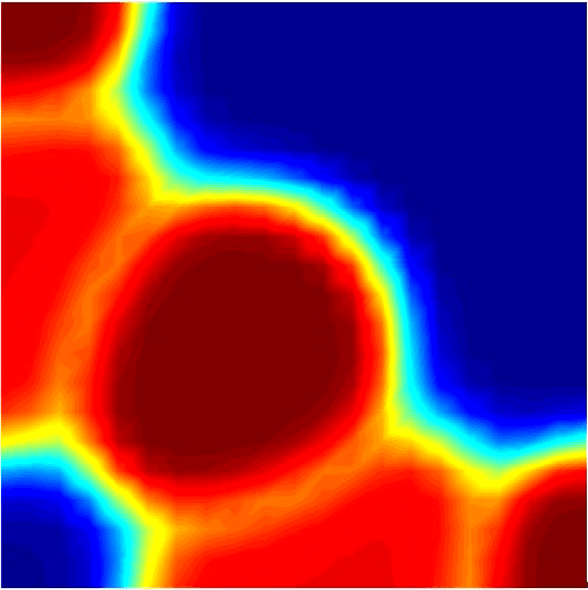}
      \includegraphics[width=1.9cm, height=1.9cm]{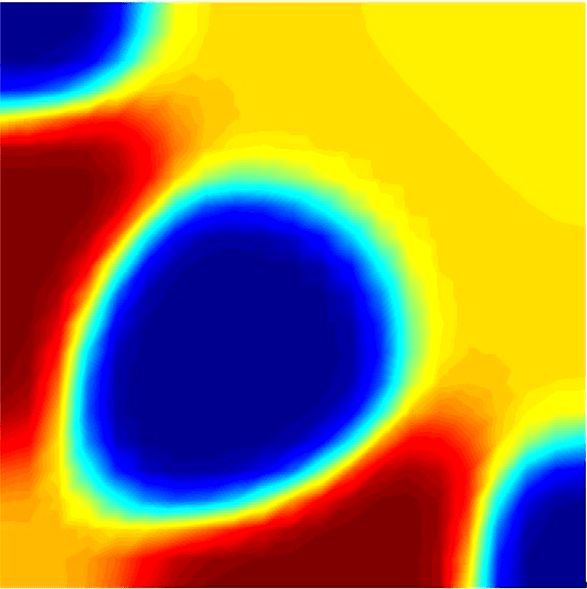}
    \end{minipage}
    \\ 
    \quad & \quad & \quad & \quad 
    \\
\hline
     \\
      \quad $50$
    & 
    \begin{minipage}{.3\textwidth}
      \includegraphics[width=1.9cm, height=1.9cm]{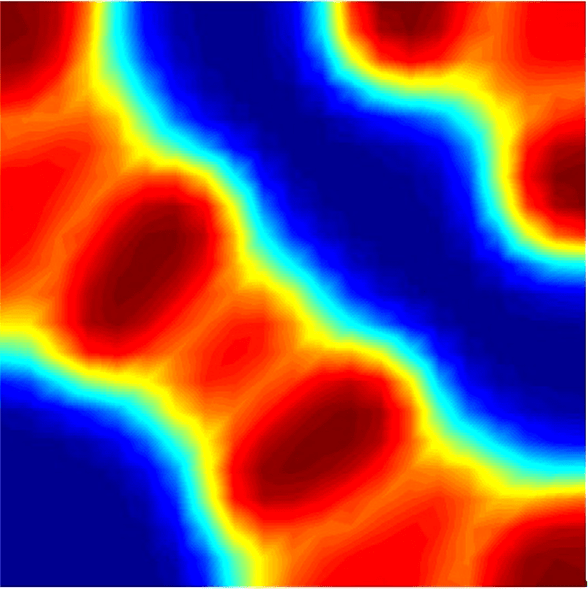}
      \includegraphics[width=1.9cm, height=1.9cm]{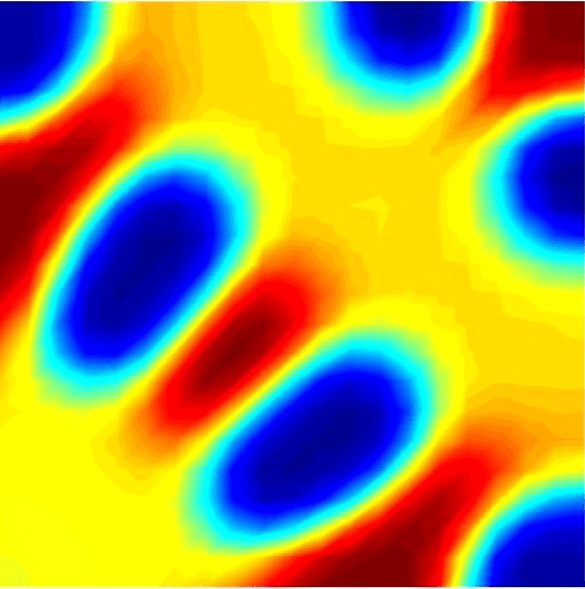}
    \end{minipage}
    &
     \begin{minipage}{.3\textwidth}
       \includegraphics[width=1.9cm, height=1.9cm]{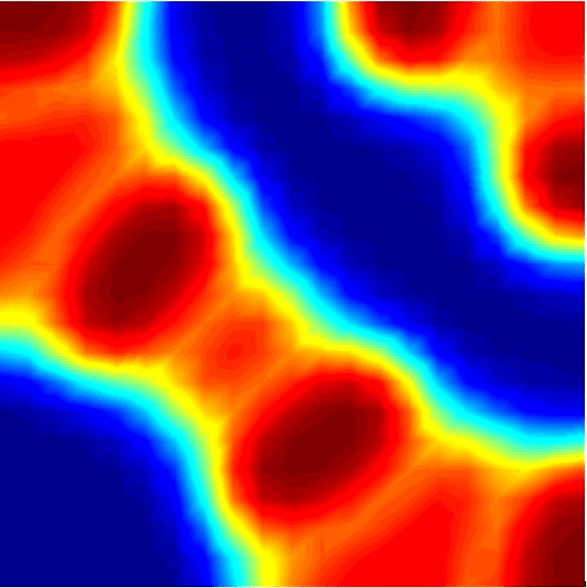}
      \includegraphics[width=1.9cm, height=1.9cm]{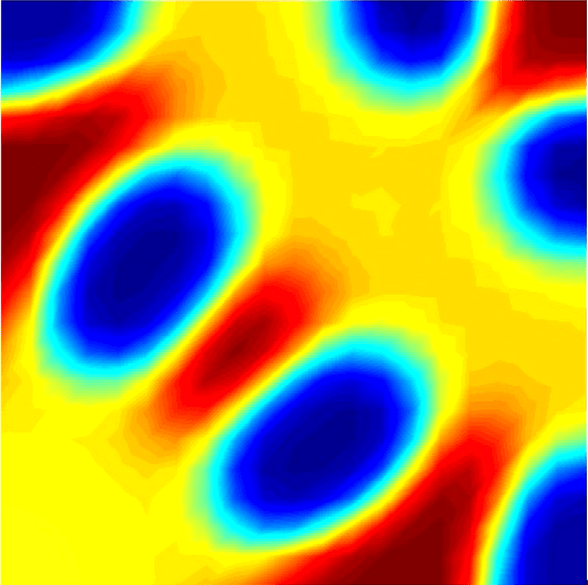}
    \end{minipage}
    &
     \begin{minipage}{.3\textwidth}
        \includegraphics[width=1.9cm, height=1.9cm]{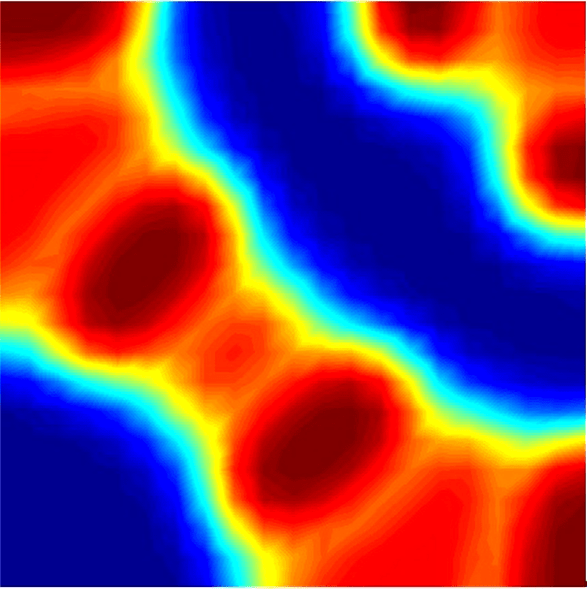}
      \includegraphics[width=1.9cm, height=1.9cm]{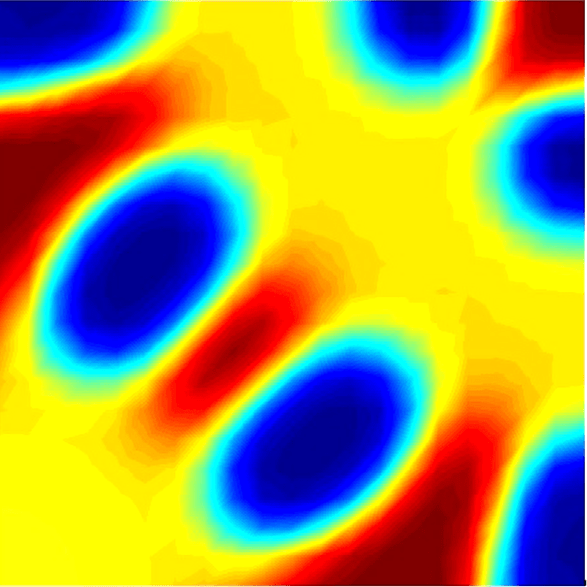}
    \end{minipage}
    \\ 
    \quad & \quad & \quad & \quad 
    \\
    \hline
     \quad & \quad & \quad & \quad 
     \\
        \quad $150$
    & 
    \begin{minipage}{.3\textwidth}
      \includegraphics[width=1.9cm, height=1.9cm]{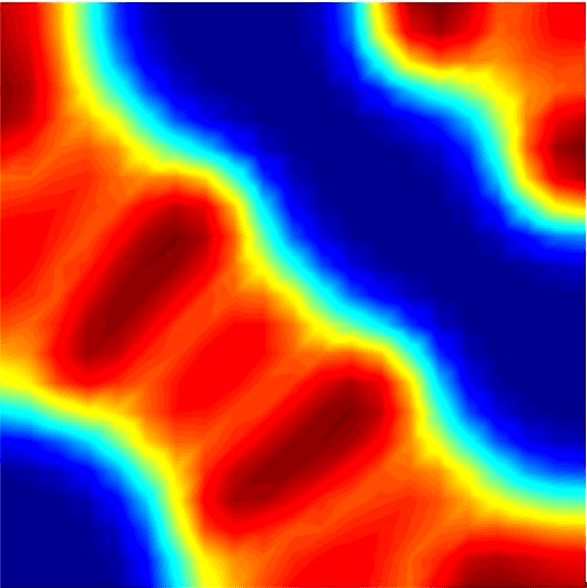}
      \includegraphics[width=1.9cm, height=1.9cm]{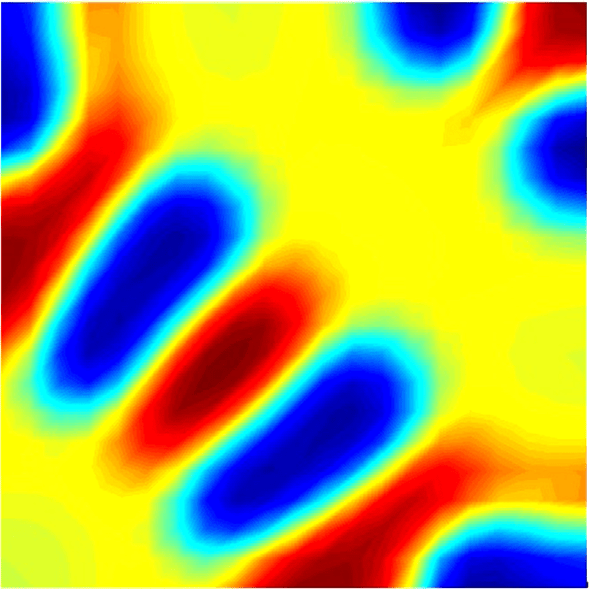}
    \end{minipage}
    &
     \begin{minipage}{.3\textwidth}
       \includegraphics[width=1.9cm, height=1.9cm]{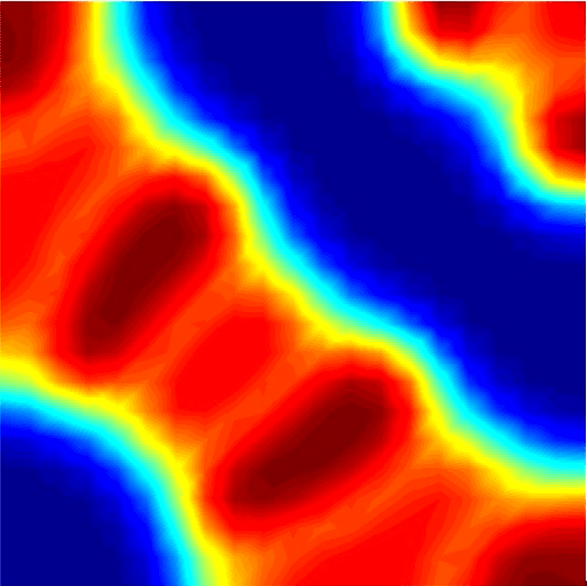}
      \includegraphics[width=1.9cm, height=1.9cm]{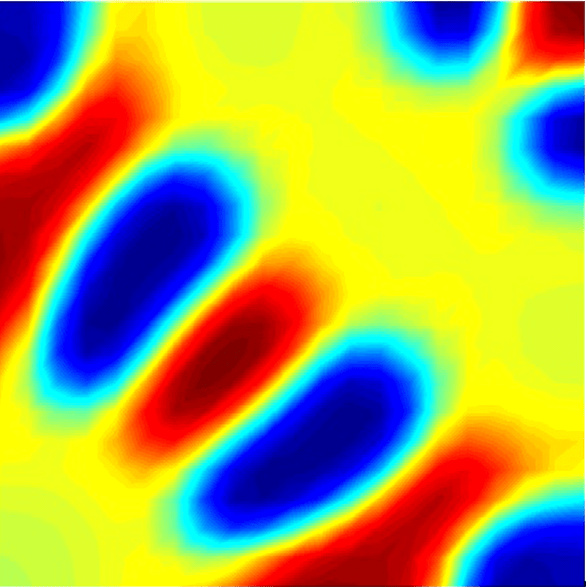}
    \end{minipage}
    &
     \begin{minipage}{.3\textwidth}
      \includegraphics[width=1.9cm, height=1.9cm]{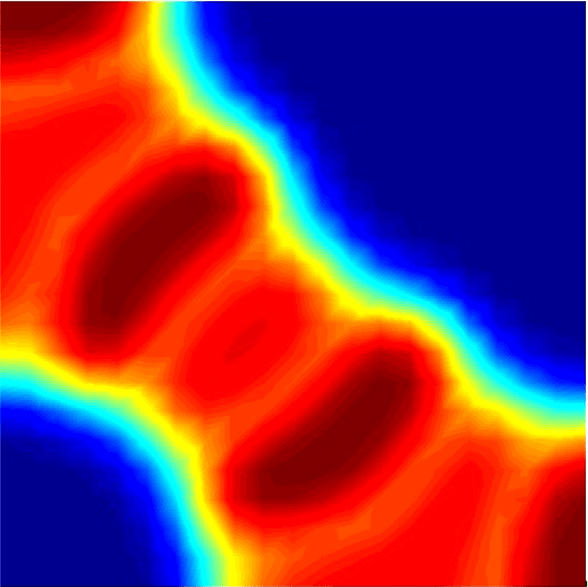}
      \includegraphics[width=1.9cm, height=1.9cm]{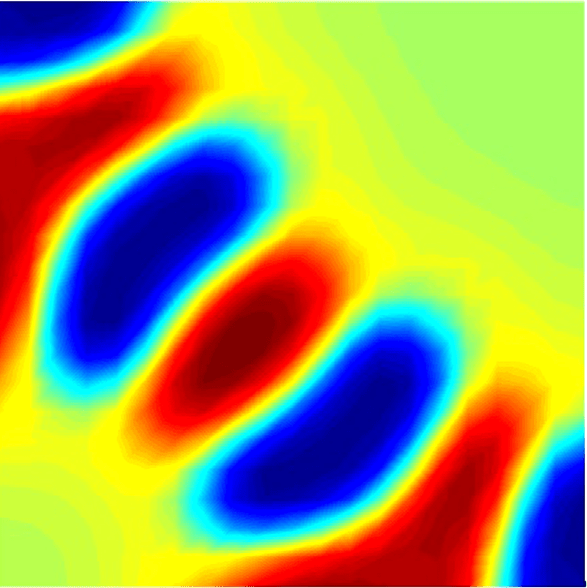}
    \end{minipage} 
    \\
    \quad & \quad & \quad & \quad \\

    \hline 
        \quad &
\quad\quad\quad\quad\quad\quad\quad\quad\quad\quad\quad 
\begin{minipage}{.3\textwidth}
       \includegraphics[trim=0 0 0 -2, scale=0.35]{colorbar.png}
    \end{minipage}
    &
    \quad 
    &
    \quad \\
  \end{tabular}
  \caption{Evolution of $u$ (left) and $v$ (right) with parameters $\epsilon_u=\epsilon_v=0.05$, $\tau_u=1$, $\tau_v=100$, $\alpha=0.4$ and $\beta=-0.9$.} 
  \label{table: sigma}
\end{table}
\quad \\
\pagebreak
\quad \\
\quad \\
\textbf{Test on $\beta$}\\
This parameter controls the energetic interaction between $u$ and $v$. As we can see in \autoref{table: beta}, if $\beta$ decreases there is a more neat separation of the two macrophases.

 \begin{table}[h!]
  \centering
  \begin{tabular}{   m{1cm}  m{5cm} m{5cm}  m{5.3cm} }
    \hline
     $\quad \beta$ &  \quad \quad \quad  \quad   $t=3$ &   \quad \quad \quad \quad $t=6$  &  \quad \quad \quad \quad  $t=15$ \\ \hline
     \quad & \quad & \quad & \quad 
     \\
      \quad $-0.3$
    & 
    \begin{minipage}{.3\textwidth}
      \includegraphics[width=1.9cm, height=1.9cm]{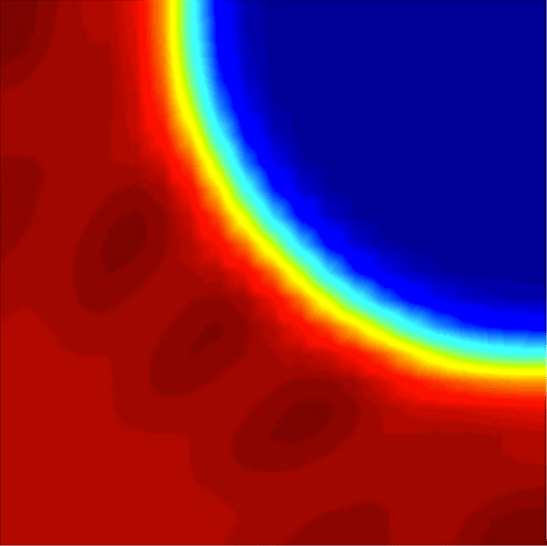}
      \includegraphics[width=1.9cm, height=1.9cm]{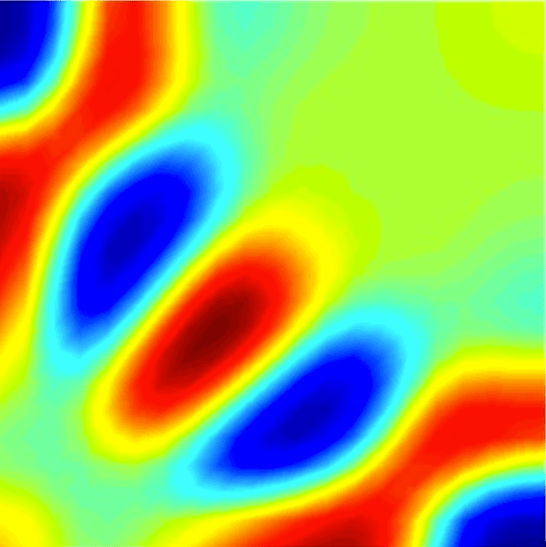}
    \end{minipage}
    &
     \begin{minipage}{.3\textwidth}
         \includegraphics[width=1.9cm, height=1.9cm]{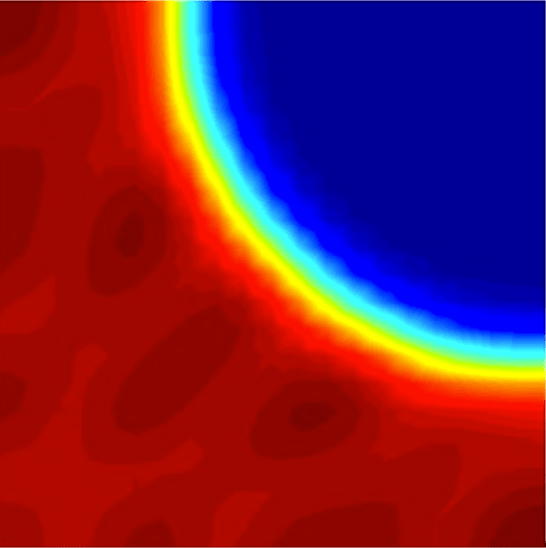}
      \includegraphics[width=1.9cm, height=1.9cm]{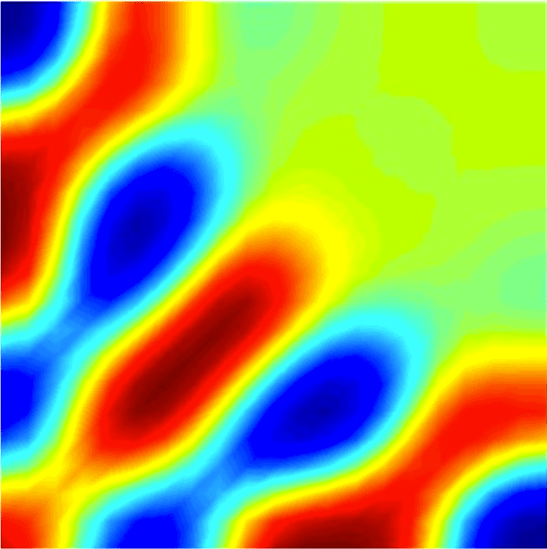}
    \end{minipage}
    &
     \begin{minipage}{.3\textwidth}
         \includegraphics[width=1.9cm, height=1.9cm]{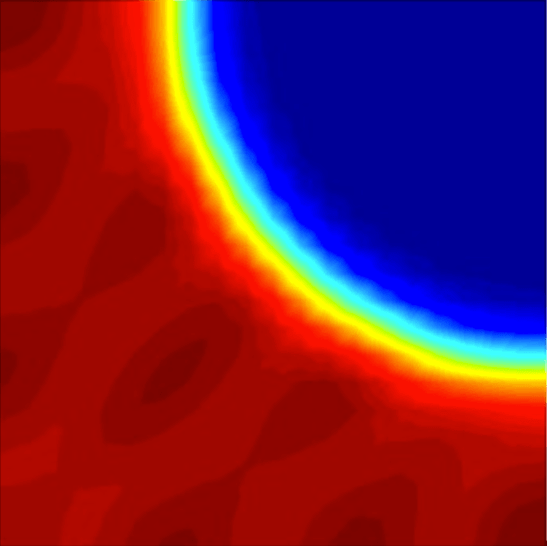}
      \includegraphics[width=1.9cm, height=1.9cm]{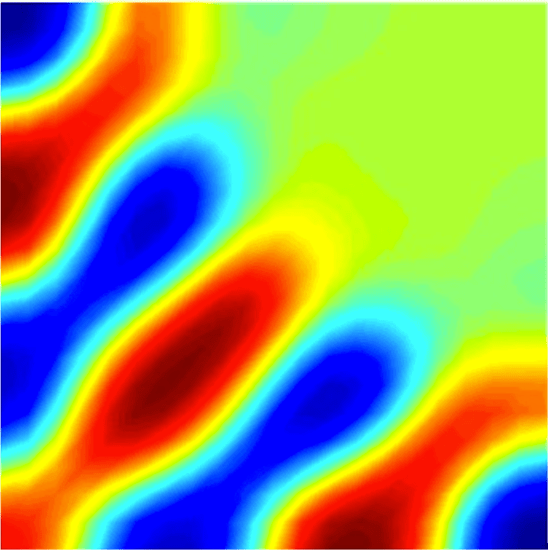}
    \end{minipage}
    \\ 
    \quad & \quad & \quad & \quad 
    \\
   
    \hline
     \quad & \quad & \quad & \quad 
    \\
    \quad $-0.5$
    & 
    \begin{minipage}{.3\textwidth}
      \includegraphics[width=1.9cm, height=1.9cm]{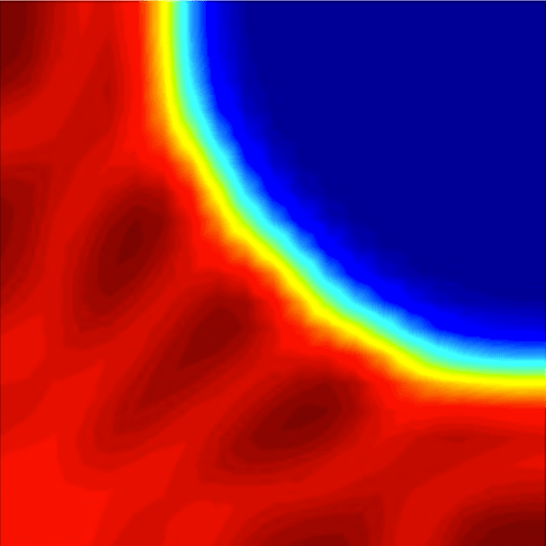}
      \includegraphics[width=1.9cm, height=1.9cm]{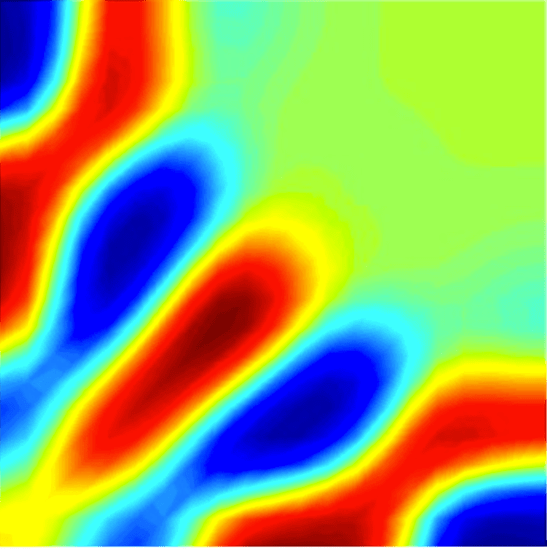}
    \end{minipage}
    &
     \begin{minipage}{.3\textwidth}
    \includegraphics[width=1.9cm, height=1.9cm]{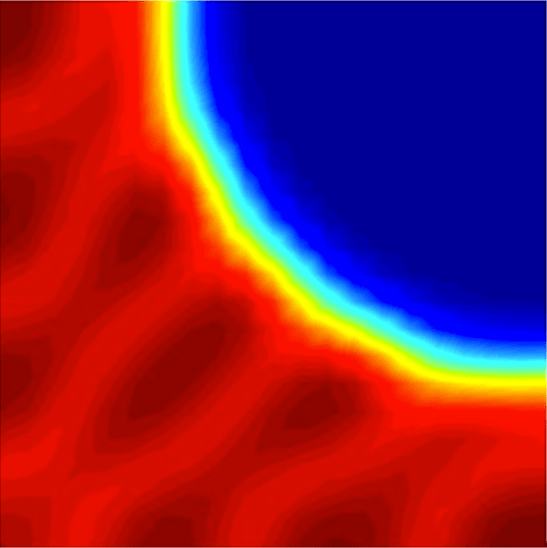}
      \includegraphics[width=1.9cm, height=1.9cm]{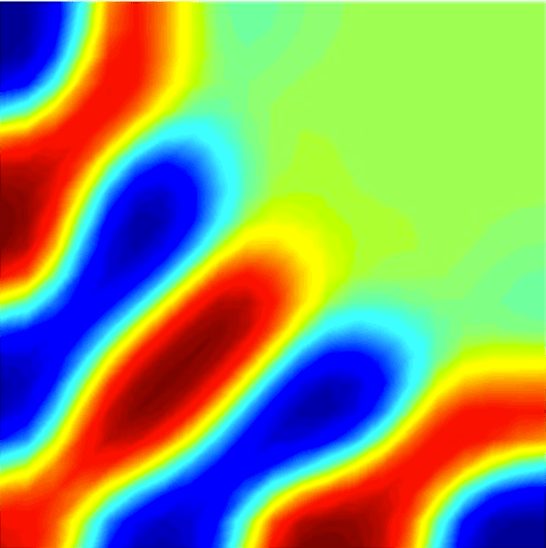}
    \end{minipage}
    &
     \begin{minipage}{.3\textwidth}
      \includegraphics[width=1.9cm, height=1.9cm]{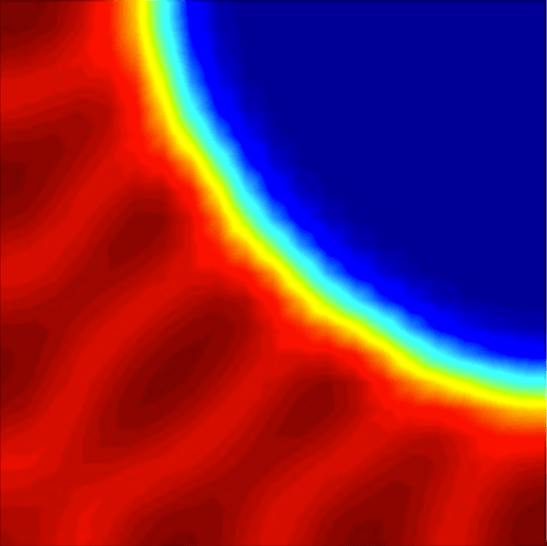}
      \includegraphics[width=1.9cm, height=1.9cm]{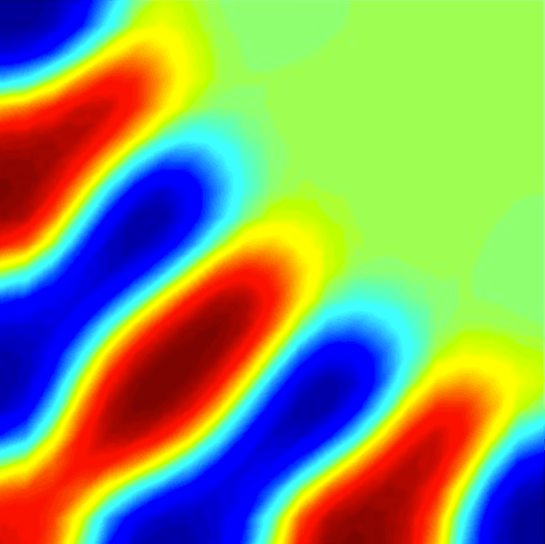}
    \end{minipage} 
    \\
    \quad & \quad & \quad & \quad \\

    \hline 
     \quad & \quad & \quad & \quad 
    \\
    \quad $-0.9$
    & 
    \begin{minipage}{.3\textwidth}
       \includegraphics[width=1.9cm, height=1.9cm]{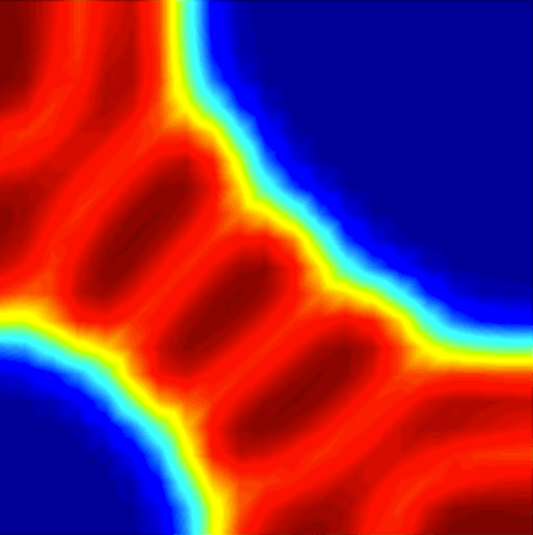}
      \includegraphics[width=1.9cm, height=1.9cm]{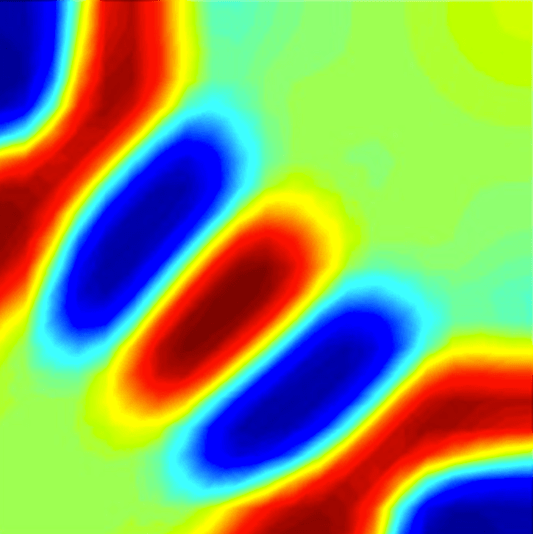}
    \end{minipage}
    &
     \begin{minipage}{.3\textwidth}
        \includegraphics[width=1.9cm, height=1.9cm]{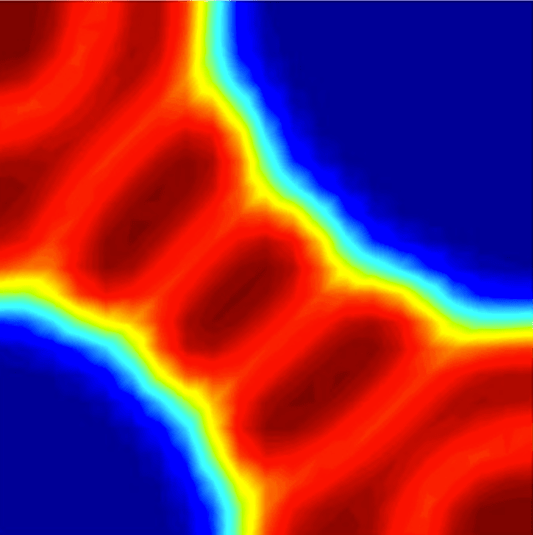}
      \includegraphics[width=1.9cm, height=1.9cm]{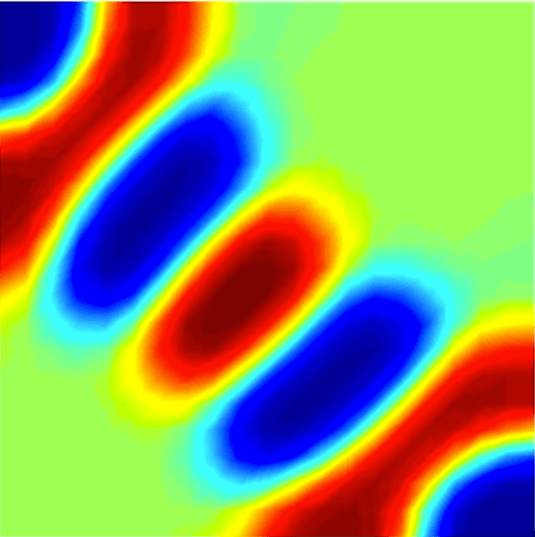}
    \end{minipage}
    &
     \begin{minipage}{.3\textwidth}
       \includegraphics[width=1.9cm, height=1.9cm]{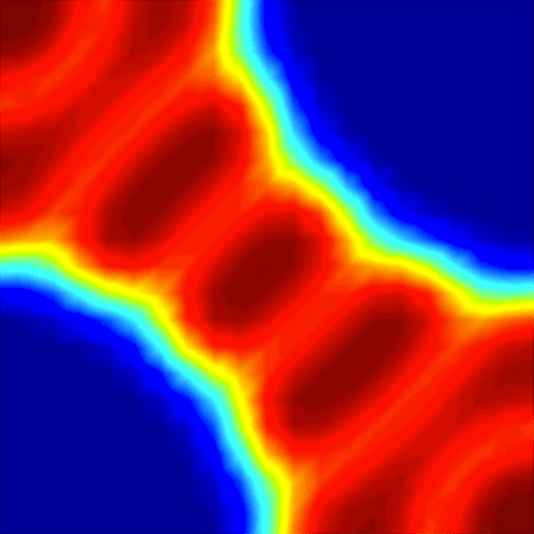}
      \includegraphics[width=1.9cm, height=1.9cm]{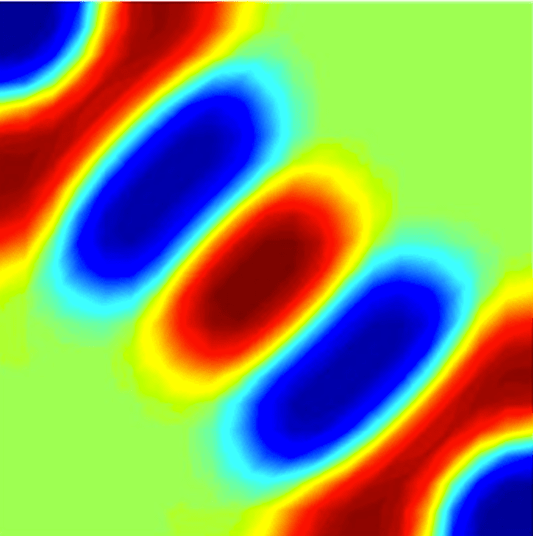}
    \end{minipage} 
    \\
    \quad & \quad & \quad & \quad \\
    
    \hline
    \quad &
\quad\quad\quad\quad\quad\quad\quad\quad\quad\quad\quad 
\begin{minipage}{.3\textwidth}
       \includegraphics[trim=0 0 0 -2, scale=0.35]{colorbar.png}
    \end{minipage}
    &
    \quad 
    &
    \quad \\
  \end{tabular}
  \caption{Evolution of $u$ (left) and $v$ (right) with parameters $\epsilon_u=\epsilon_v=0.05$, $\tau_u=1$, $\tau_v=100$, $\sigma=100$, and $\alpha=0.01$.} 
  \label{table: beta}
\end{table}
\pagebreak
\quad \\
\textbf{Test on $ \alpha$} \\
This value controls the interaction between the confined copolymer and the confining surface, causing symmetry-breaking between microphase separated domains. If $\alpha = 0$, both microphases in $v$ can reach $u$ with the same probability, but if $\alpha \neq 0$ we are changing the preference of $u$ for the positive or negative values of $v$ and thus the confining surface will change according to this preference. Indeed, as we can see in \autoref{table: alfa}, for negative values of $\alpha$, at the interface with $u$ we have prevalence of $v=-1$, instead, with $\alpha$ positive, at the interface with $u$ we have prevalence of $v=1$.

 \begin{table}[h!]
  \centering
  \begin{tabular}{   m{1cm}  m{5cm} m{5cm}  m{5.3cm} }
    \hline
     $\quad \alpha$ &  \quad \quad \quad  \quad   $t=3$ &  \quad \quad \quad \quad  $t=6$  &   \quad \quad \quad \quad  $t=15$ \\ \hline
     \quad & \quad & \quad & \quad 
     \\
      \quad $-0.4$
    & 
    \begin{minipage}{.3\textwidth}
      \includegraphics[width=1.9cm, height=1.9cm]{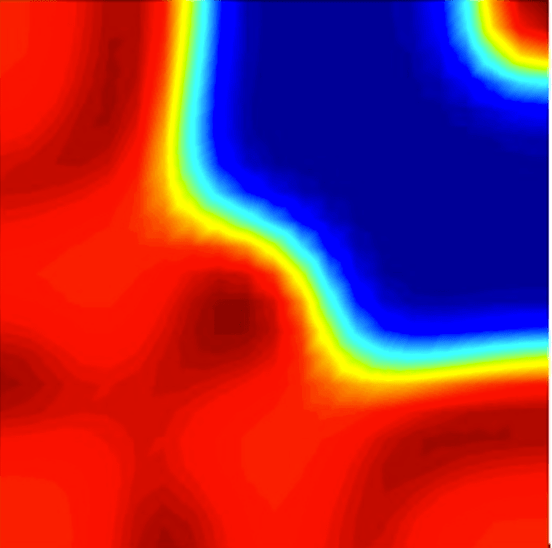}
      \includegraphics[width=1.9cm, height=1.9cm]{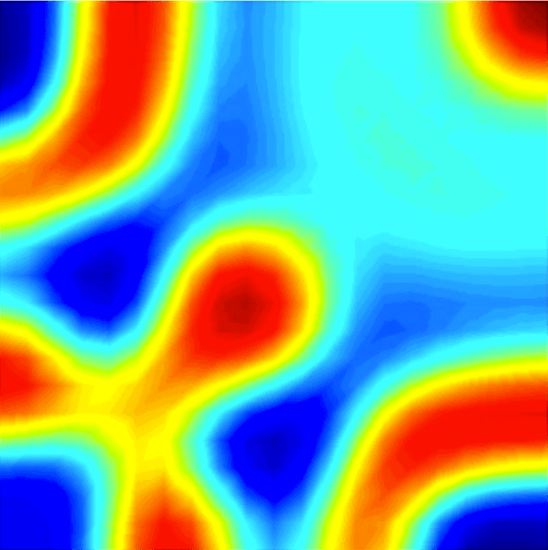}
    \end{minipage}
    &
     \begin{minipage}{.3\textwidth}
       \includegraphics[width=1.9cm, height=1.9cm]{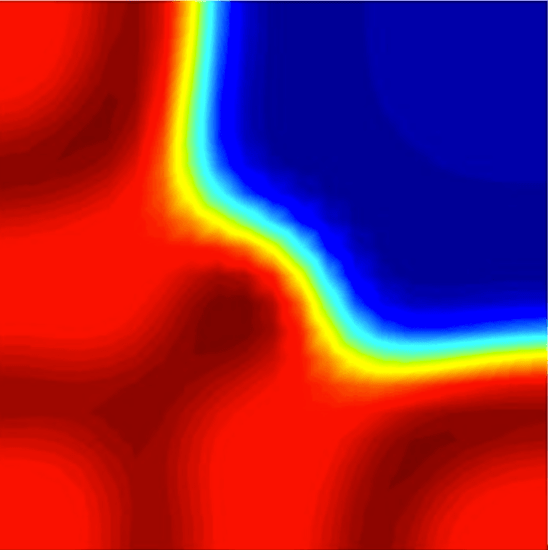}
      \includegraphics[width=1.9cm, height=1.9cm]{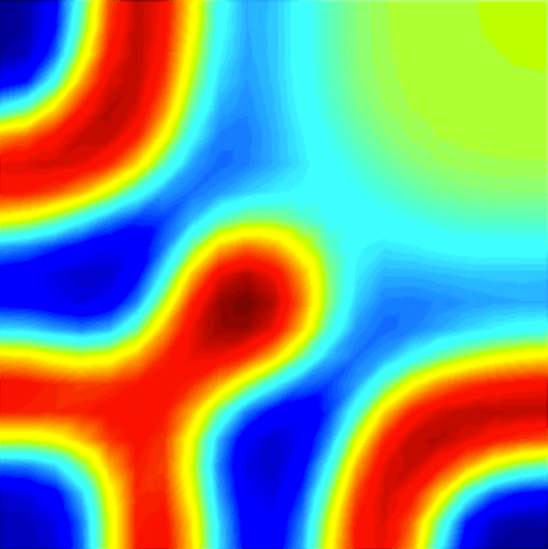}
    \end{minipage}
    &
     \begin{minipage}{.3\textwidth}
       \includegraphics[width=1.9cm, height=1.9cm]{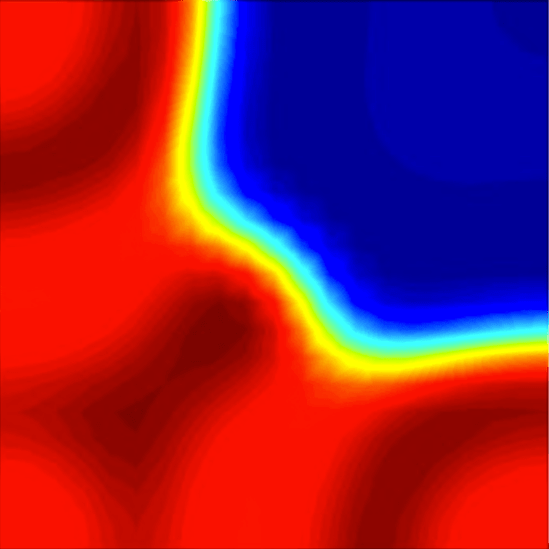}
      \includegraphics[width=1.9cm, height=1.9cm]{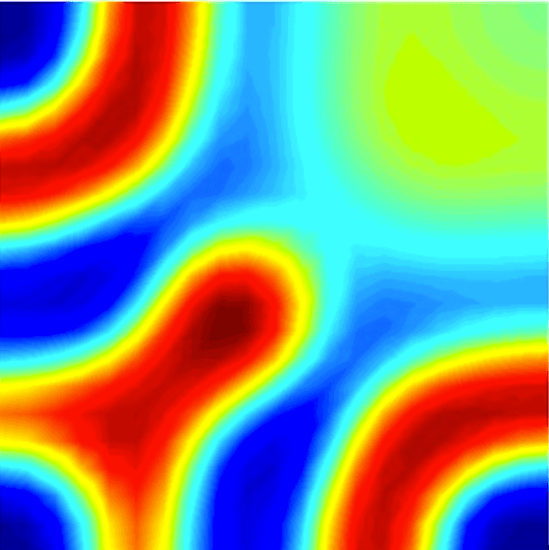}
    \end{minipage}\\
    \quad & \quad & \quad & \quad \\
    \hline
     \quad & \quad & \quad & \quad \\ \quad $-0.08$
    & 
    \begin{minipage}{.3\textwidth}
      \includegraphics[width=1.9cm, height=1.9cm]{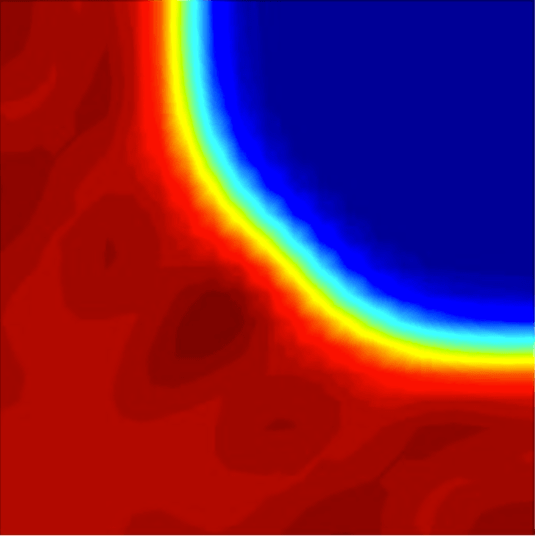}
      \includegraphics[width=1.9cm, height=1.9cm]{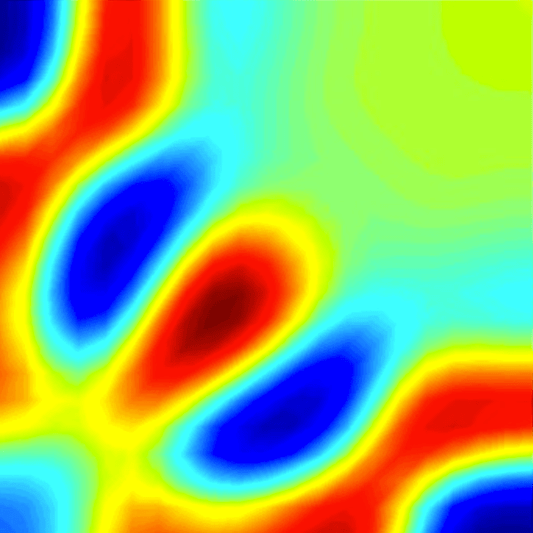}
    \end{minipage}
    &
     \begin{minipage}{.3\textwidth}
       \includegraphics[width=1.9cm, height=1.9cm]{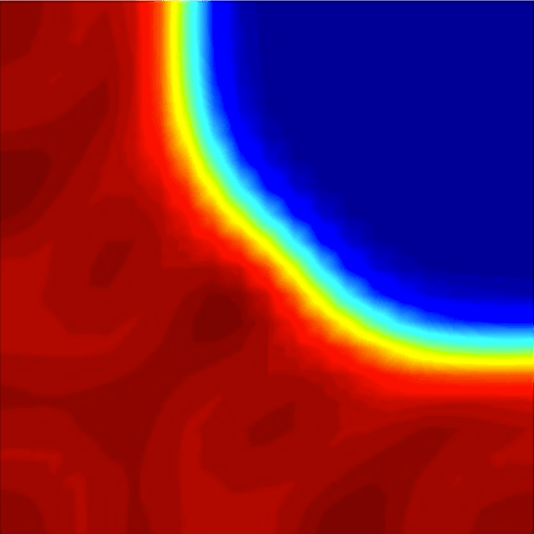}
      \includegraphics[width=1.9cm, height=1.9cm]{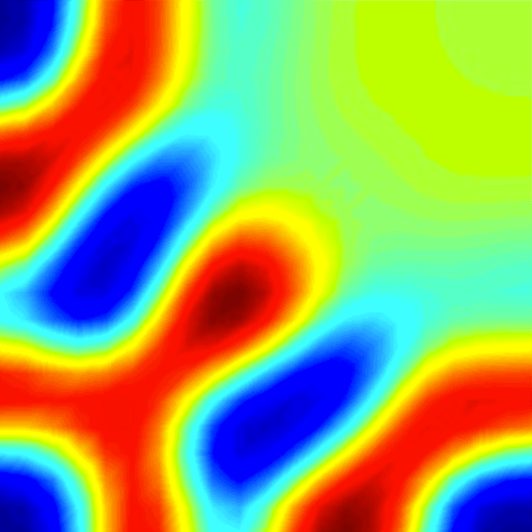}
    \end{minipage}
    &
     \begin{minipage}{.3\textwidth}
       \includegraphics[width=1.9cm, height=1.9cm]{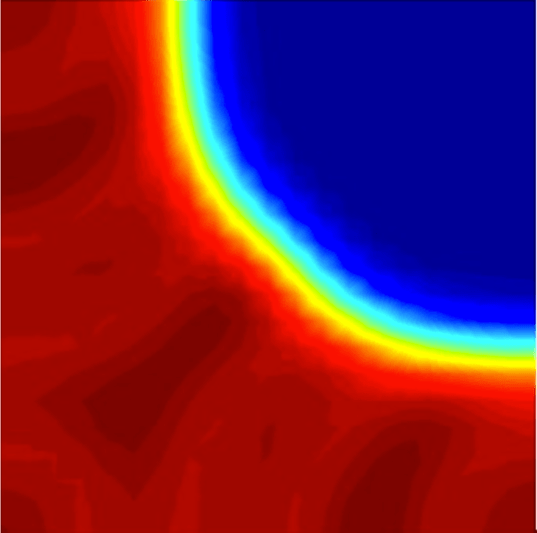}
      \includegraphics[width=1.9cm, height=1.9cm]{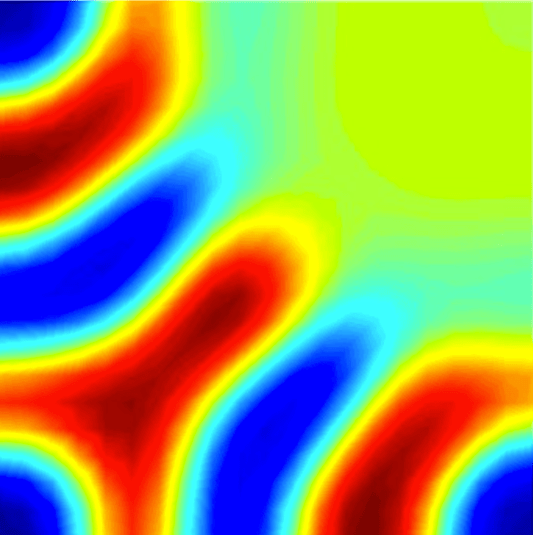}
    \end{minipage}\\
    \quad & \quad & \quad & \quad \\
    \hline
     \quad & \quad & \quad & \quad \\
      \quad $0.01$
    & 
    \begin{minipage}{.3\textwidth}
      \includegraphics[width=1.9cm, height=1.9cm]{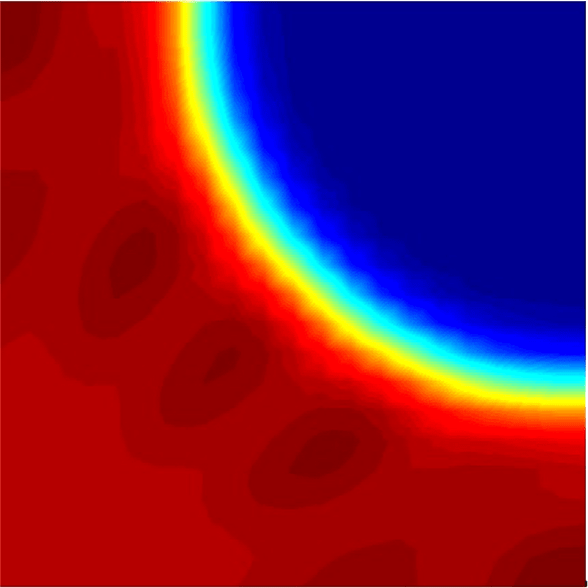}
      \includegraphics[width=1.9cm, height=1.9cm]{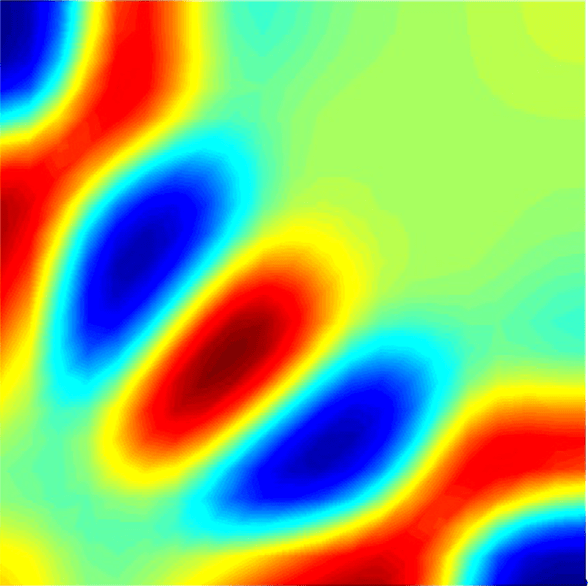}
    \end{minipage}
    &
     \begin{minipage}{.3\textwidth}
       \includegraphics[width=1.9cm, height=1.9cm]{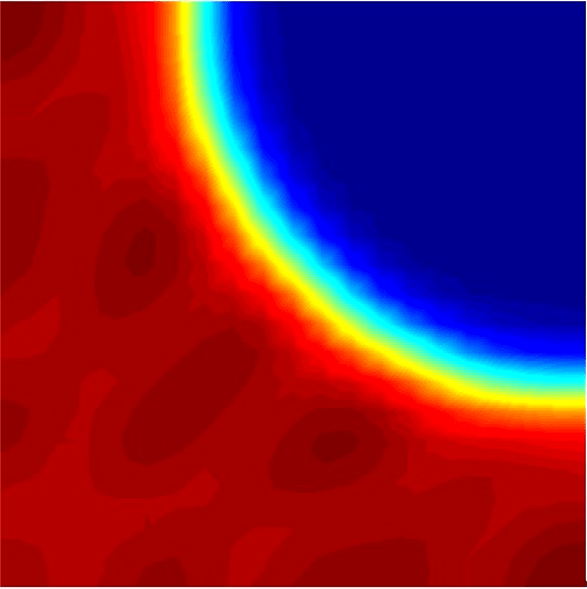}
      \includegraphics[width=1.9cm, height=1.9cm]{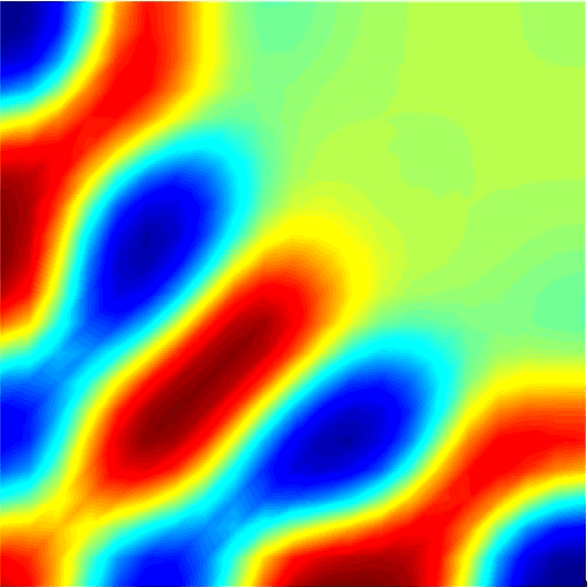}
    \end{minipage}
    &
     \begin{minipage}{.3\textwidth}
       \includegraphics[width=1.9cm, height=1.9cm]{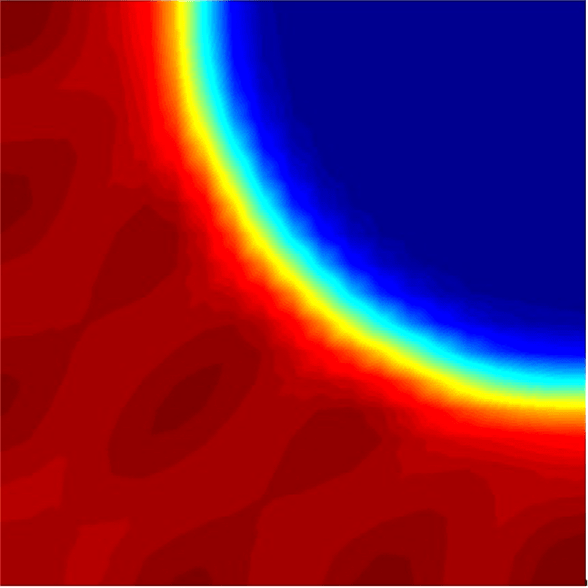}
      \includegraphics[width=1.9cm, height=1.9cm]{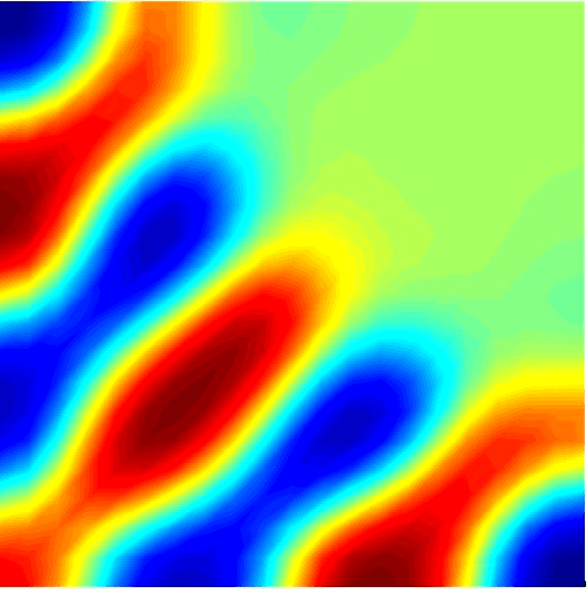}
    \end{minipage}\\
    \quad & \quad & \quad & \quad \\
    \hline
     \quad & \quad & \quad & \quad \\
\quad $0.1$
    & 
    \begin{minipage}{.3\textwidth}
      \includegraphics[width=1.9cm, height=1.9cm]{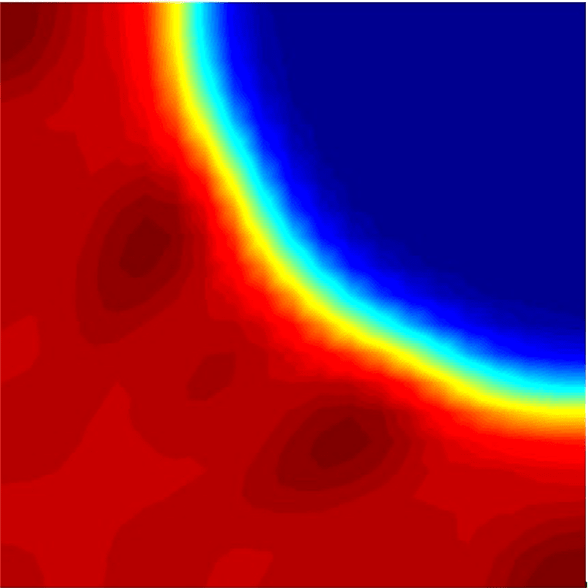}
      \includegraphics[width=1.9cm, height=1.9cm]{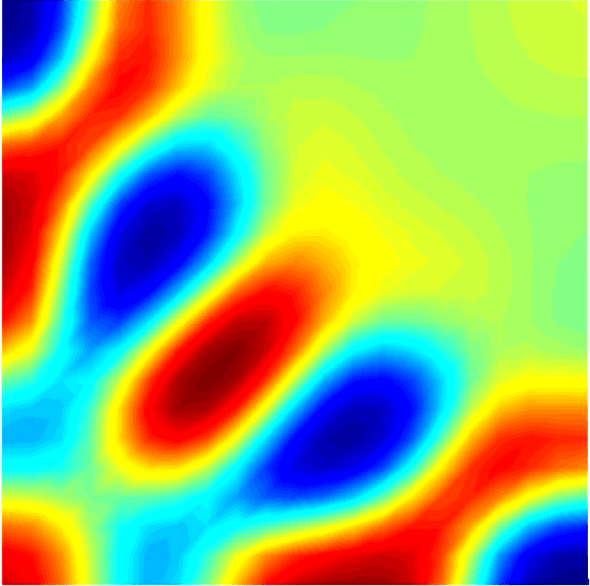}
    \end{minipage}
    &
     \begin{minipage}{.3\textwidth}
       \includegraphics[width=1.9cm, height=1.9cm]{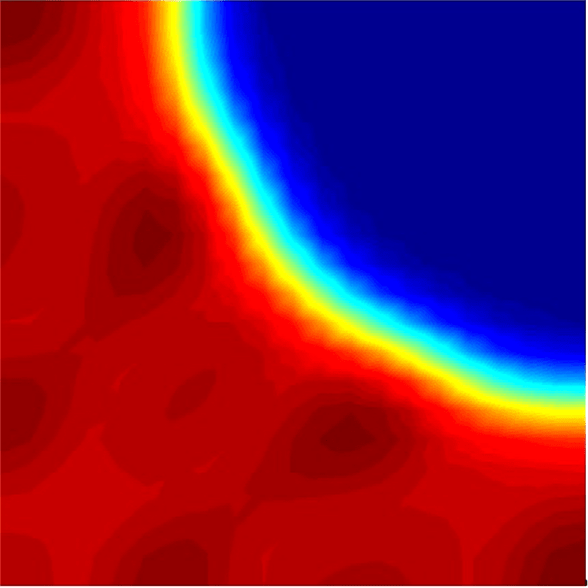}
      \includegraphics[width=1.9cm, height=1.9cm]{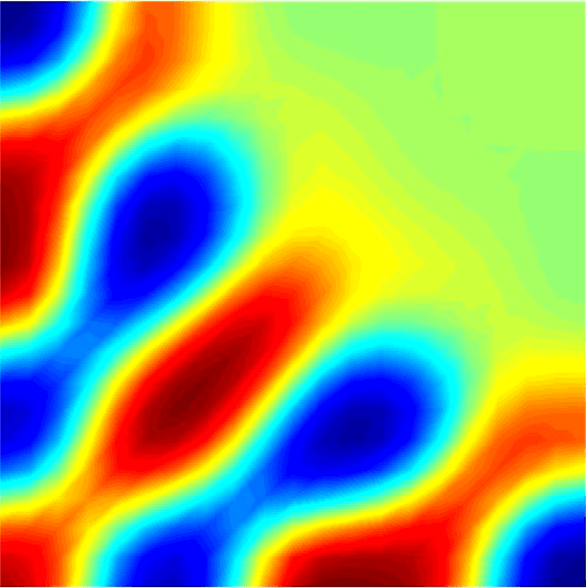}
    \end{minipage}
    &
     \begin{minipage}{.3\textwidth}
       \includegraphics[width=1.9cm, height=1.9cm]{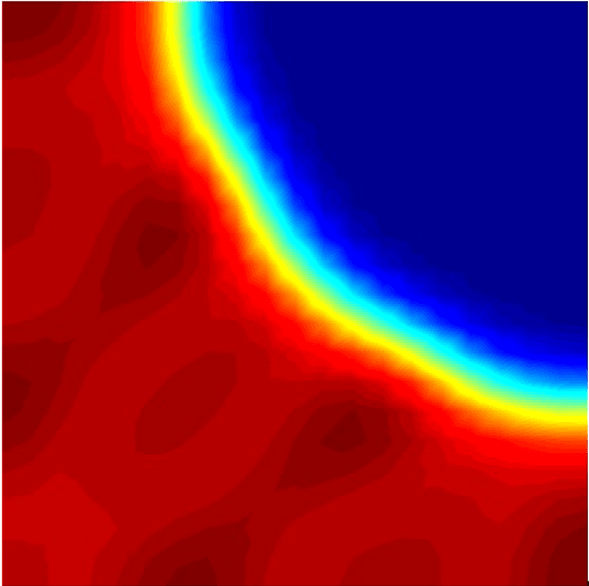}
      \includegraphics[width=1.9cm, height=1.9cm]{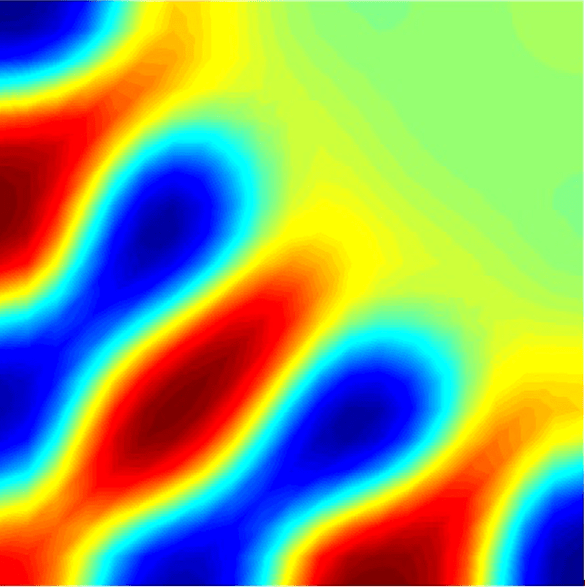}
    \end{minipage}\\
    \quad & \quad & \quad & \quad \\
    \hline
     \quad & \quad & \quad & \quad \\    \quad $0.4$
    & 
   \begin{minipage}{.3\textwidth}
      \includegraphics[width=1.9cm, height=1.9cm]{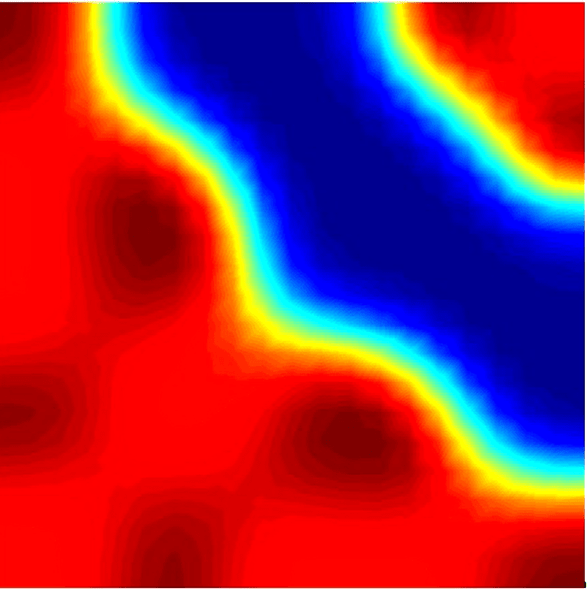}
      \includegraphics[width=1.9cm, height=1.9cm]{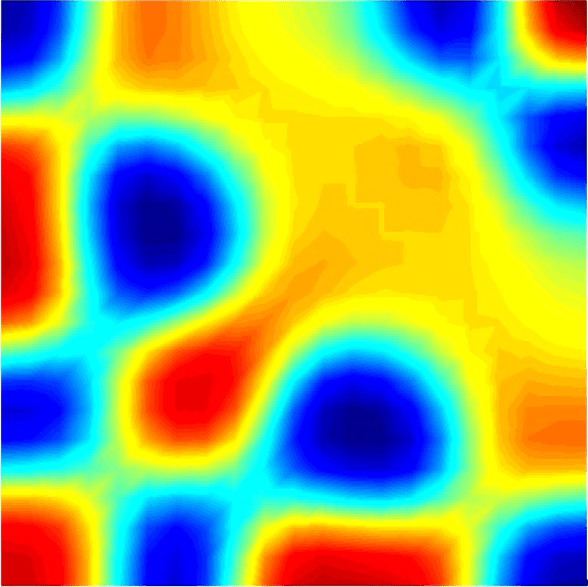}
    \end{minipage}
    &
     \begin{minipage}{.3\textwidth}
      \includegraphics[width=1.9cm, height=1.9cm]{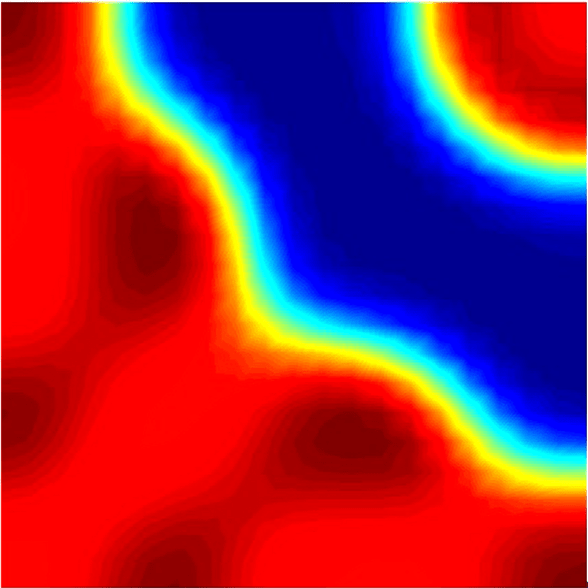}
      \includegraphics[width=1.9cm, height=1.9cm]{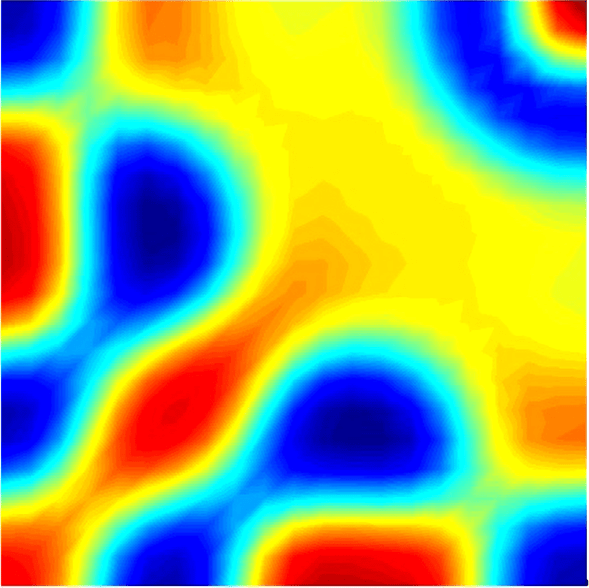}
    \end{minipage}
    &
     \begin{minipage}{.3\textwidth}
      \includegraphics[width=1.9cm, height=1.9cm]{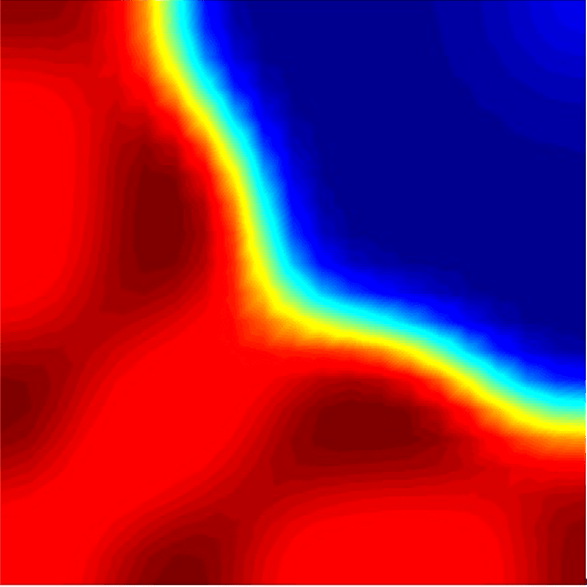}
      \includegraphics[width=1.9cm, height=1.9cm]{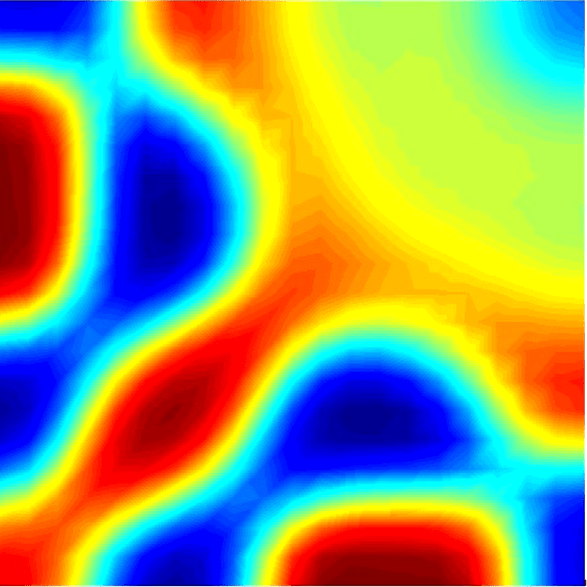}
    \end{minipage} \\
    \quad & \quad & \quad & \quad \\
    \hline
        \quad &
\quad\quad\quad\quad\quad\quad\quad\quad\quad\quad\quad 
\begin{minipage}{.3\textwidth}
       \includegraphics[trim=0 0 0 -2, scale=0.35]{colorbar.png}
    \end{minipage}
    &
    \quad 
    &
    \quad \\

  \end{tabular}
  \caption{Evolution of $u$ (left) and $v$ (right) with parameters $\epsilon_u=\epsilon_v=0.05$, $\tau_u=1$, $\tau_v=100$, $\sigma=100$, and $\beta=-0.3$.}
  \label{table: alfa}
\end{table}
\pagebreak
\quad \\
\textbf{Mass conservation} \\
Finally, we verify that if $ \overline{v}= \int_{\Omega}v_0 $, the quantity $ \int_{\Omega} v(t)$ remains constant during the evolution. On the contrary, if we set an arbitrary value of $ \overline{v}$ we can see that there is no mass conservation but $ \int_{\Omega} v(t) \to \overline{v} $ as $ t \to + \infty$. We report the following test in \autoref{table: vsegnato} with setting: $\Omega= (0,1) \times (0,1)$, $u_0= \sin(10xy)$, $v_0= \cos(10(x-y))xy$, $\Delta t=0.005$, $T=15$ and the FEM space we use is still $P1$ on a $20\times20$ grid. 

  \begin{table}[h!]
  \centering
  \begin{tabular}{   m{1cm}  m{5cm} m{5cm}  m{5.3cm} }
    \hline
     $\quad \overline{v}$ &  \quad \quad  \quad   $t=1$ &   \quad \quad \quad $t=7.5$  &  \quad \quad \quad  $t=15$ \\ \hline
     \quad & \quad & \quad & \quad 
     \\
      \quad $0$
    & 
    \begin{minipage}{.3\textwidth}
      \includegraphics[width=3cm, height=3cm]{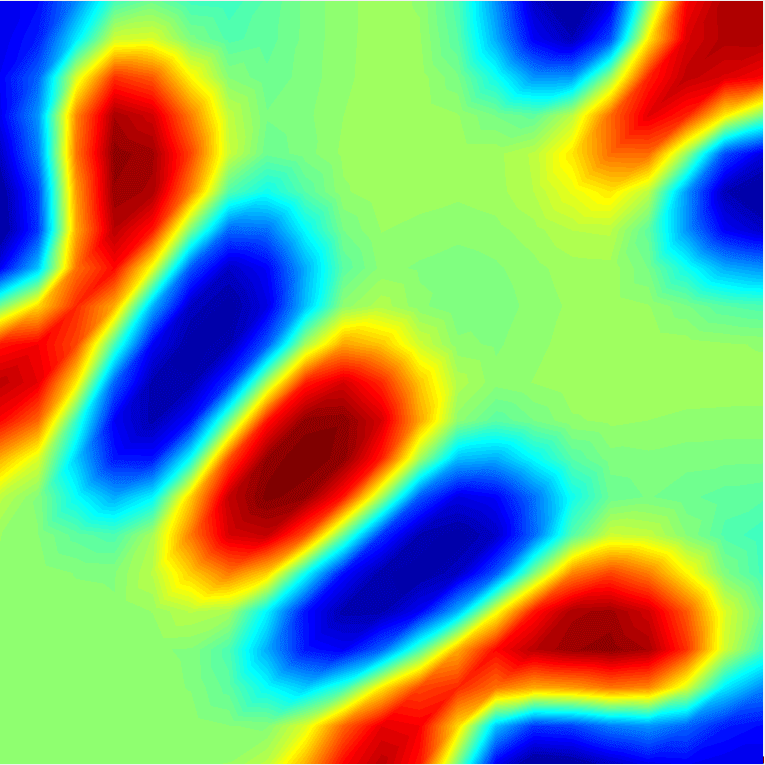}
      \\ $ \int_{\Omega} v(t) = 0.0042$
    \end{minipage}
    &
     \begin{minipage}{.3\textwidth}
         \includegraphics[width=3cm, height=3cm]{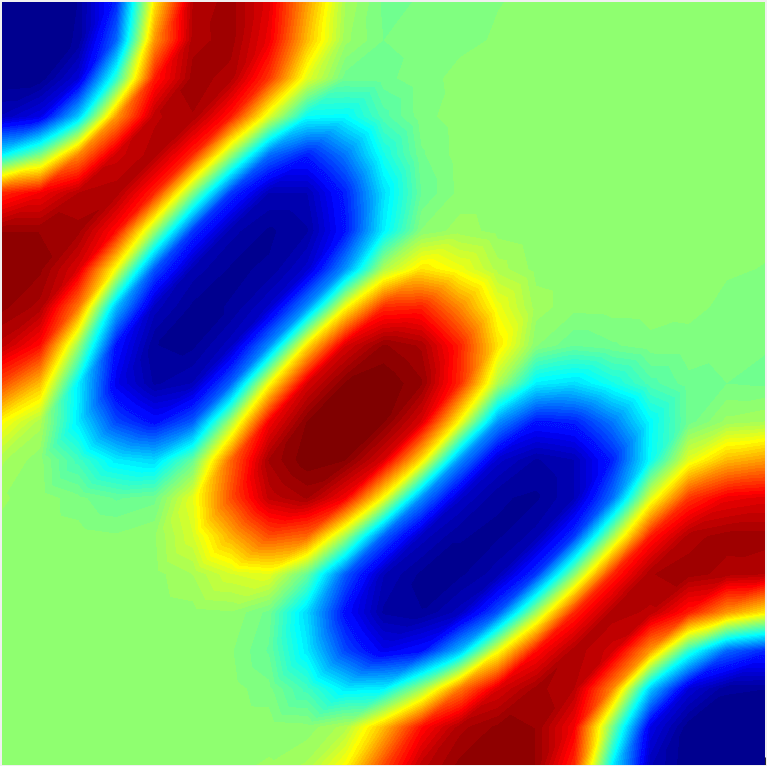}
         \\ $ \int_{\Omega} v(t) \approx 10^{-6}$
         \end{minipage}
    &
     \begin{minipage}{.3\textwidth}
          \includegraphics[width=3cm, height=3cm]{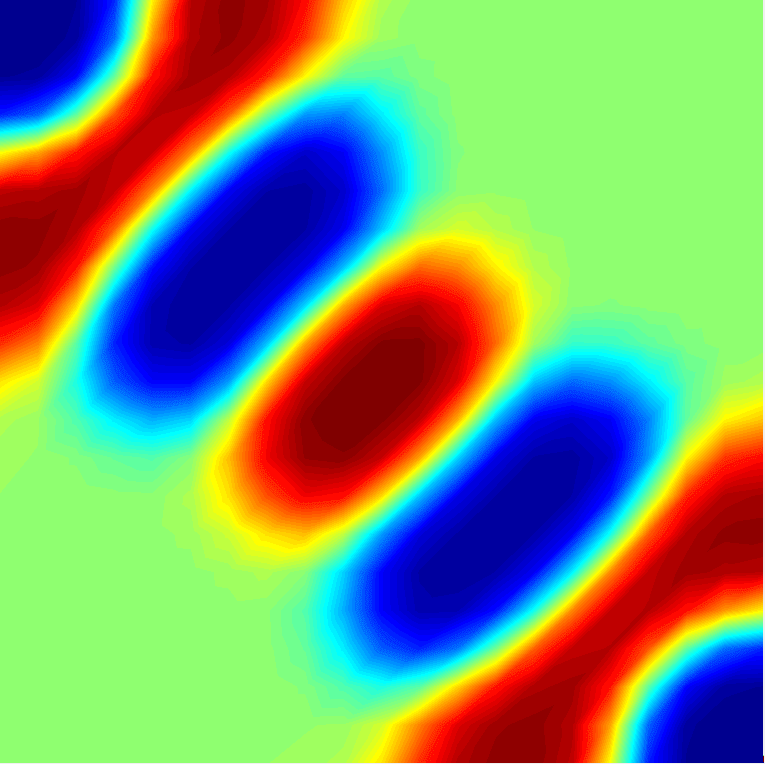}
          \\ $ \int_{\Omega} v(t) \approx 10^{-7} $
    \end{minipage}
    \\ 
    \quad & \quad & \quad & \quad 
    \\
   
    \hline
     \quad & \quad & \quad & \quad 
    \\
    \quad $\int_{\Omega}v_0$
    & 
    \begin{minipage}{.3\textwidth}
      \includegraphics[width=3cm, height=3cm]{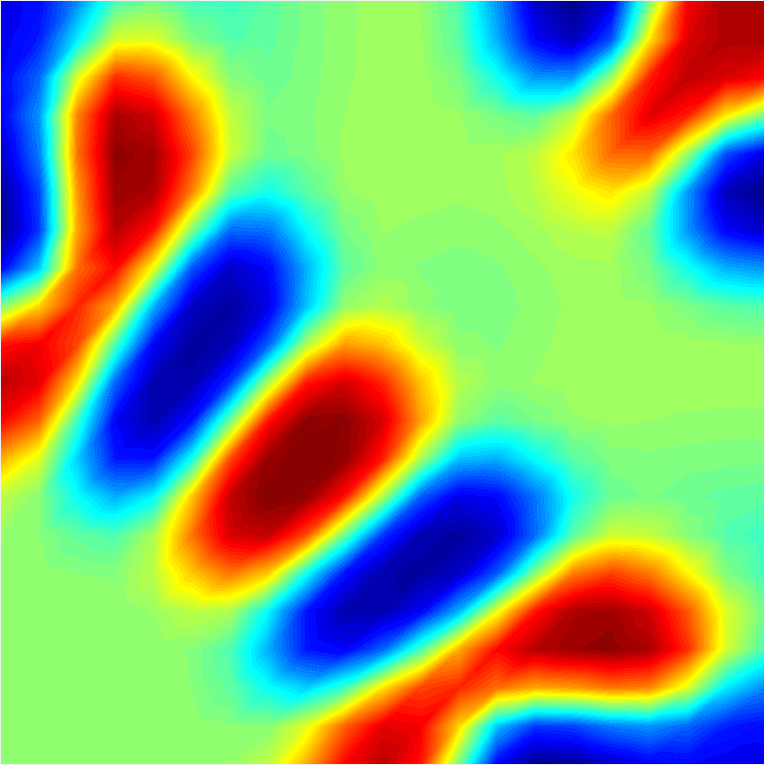}
      \\ $ \int_{\Omega} v(t) = 0.01139$
      
    \end{minipage}
    &
     \begin{minipage}{.3\textwidth}
    \includegraphics[width=3cm, height=3cm]{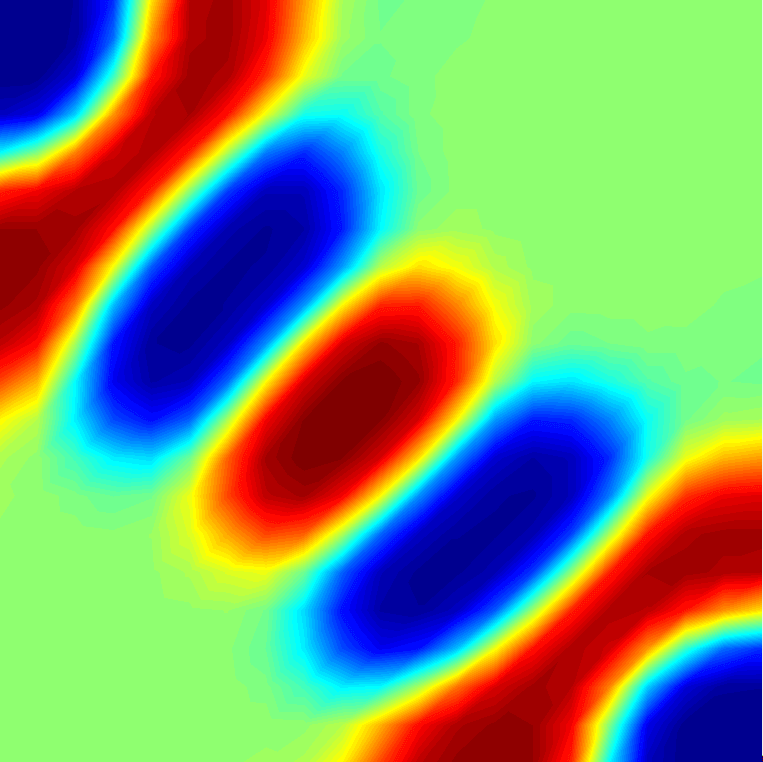}
    \\ $ \int_{\Omega} v(t) =0.01146$
    \end{minipage}
    &
     \begin{minipage}{.3\textwidth}
      \includegraphics[width=3cm, height=3cm]{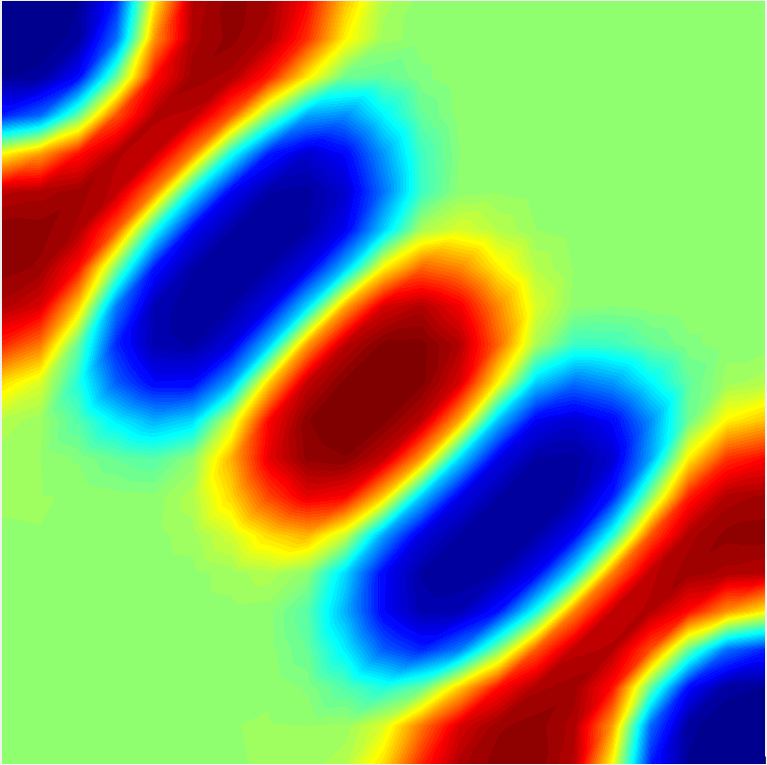}
      \\ $ \int_{\Omega} v(t) =0.01146 $
    \end{minipage} 
    \\
    \quad & \quad & \quad & \quad \\

    \hline 
     \quad & \quad & \quad & \quad 
    \\
    \quad $0.6$
    & 
    \begin{minipage}{.3\textwidth}
       \includegraphics[width=3cm, height=3cm]{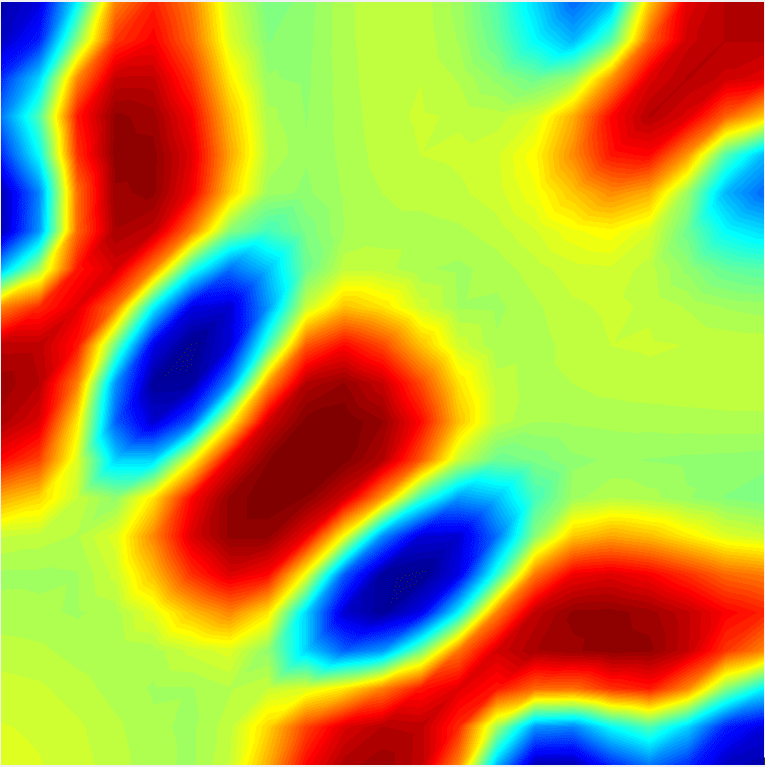}
       \\ $ \int_{\Omega} v(t) =0.38179 $
    \end{minipage}
    &
     \begin{minipage}{.3\textwidth}
        \includegraphics[width=3cm, height=3cm]{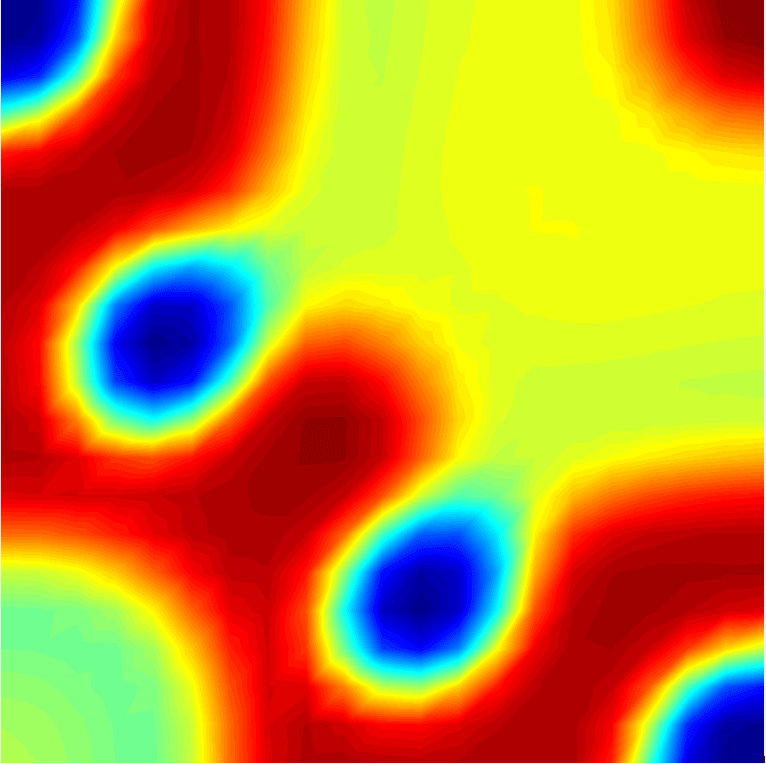}
        \\ $ \int_{\Omega} v(t) =0.59967 $
    \end{minipage}
    &
     \begin{minipage}{.3\textwidth}
       \includegraphics[width=3cm, height=3cm]{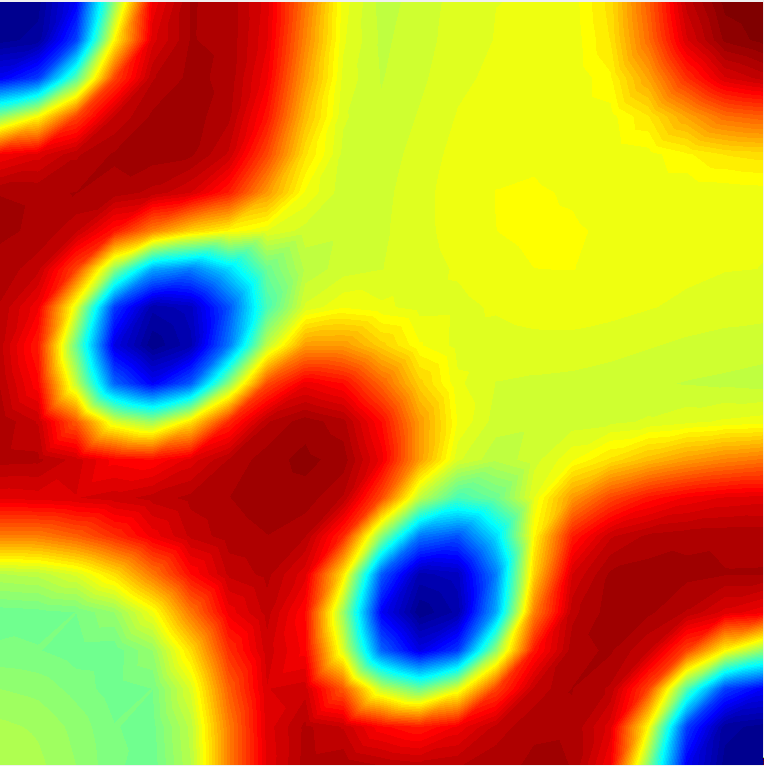}
       \\ $ \int_{\Omega} v(t) =0.59999 $
    \end{minipage} 
    \\
    \quad & \quad & \quad & \quad \\
    
    \hline
    \quad &
\quad\quad\quad\quad\quad\quad\quad\quad\quad\quad
\begin{minipage}{.3\textwidth}
       \includegraphics[trim=0 0 0 -2, scale=0.35]{colorbar.png}
    \end{minipage}
    &
    \quad 
    &
    \quad \\
  \end{tabular}
  \caption{Evolution of $v$ with parameters $\epsilon_u=\epsilon_v=0.05$, $\tau_u=1$, $\tau_v=100$, $\sigma=100$, and $\alpha=0.01$, $\beta=-0.9$. Notice that $\int_{\Omega}v(t)=0.01146$.} 
  \label{table: vsegnato}
\end{table}
\pagebreak
 \section{Conclusions} \label{concl}
We have proposed and tested different strategies to solve numerically a system of fully coupled Cahn-Hilliard equations. The numerical solutions we have obtained seem to give a fair description of the physical process. 
From the sensitivity analysis on the robustness with respect to the physical parameters, we can see that the evolution of the copolymer is confined where $u=1$. This is very important, since the solution behaves coherently with the phenomenon without imposing conditions at the interface.\\
Our results in the two-dimensional framework are in accordance with the ones in \cite{grass}, even if we consider Neumann boundary conditions instead of periodic. 
From a numerical point of view, we notice that the linearization methods used for the single Cahn-Hilliard equation can be used also in the context of the coupled system, with the exception of the Wu-Van Zwieten-Van der Zee's method. Moreover, the results are good even with a coarse grid in space ($20\times20$) and with a low order (i.e $1$) FEM approximations.\\
This work represents a first step in developing efficient schemes for fully coupled systems of Cahn-Hilliard equations. It may be interesting to further carry out some experimental validations and to fully address the theoretical analysis. Specifically, a further numerical development could be the theoretical analysis in the coupled case of the methods we applied. From an analytical point of view, the study of the well posedness, regularity and longtime behavior of the solutions is also a topic of further research.
\pagebreak

\printbibliography
\vfill
\end{document}